\documentclass[autocontact]{gaceta}

\usepackage[T1]{fontenc}
\usepackage{lmodern}
\usepackage[utf8]{inputenc}
\usepackage[spanish]{babel}

\let\oldmarkright\markright
\let\oldmarkboth\markboth
\renewcommand
\renewcommand

\usepackage{bbm}
\usepackage{xurl} 
\usepackage{wrapfig}

\usepackage[pdfborder={0 0 0}]{hyperref}

\newcommand{\Z}{\mathbb{Z}}
\newcommand{\R}{\mathbb{R}}

\newtheorem{definition}{Definición}
\newtheorem{lemma}{Lema}
\newtheorem{theorem}{Teorema}
\newtheorem{corollary}{Corolario}
\newtheorem{proposition}{Proposición}
\newtheorem{example}{Ejemplo}
\newtheorem{conjecture}{Conjetura}

\DeclareMathOperator{\sgn}{sgn}
\DeclareMathOperator{\Cov}{\mathbf{Cov}}
\DeclareMathOperator{\Var}{\mathbf{Var}}
\DeclareMathOperator{\II}{\mathbf{I}}
\DeclareMathOperator{\EE}{\mathbf{E}}
\DeclareMathOperator{\HH}{\mathbf{H}}
\DeclareMathOperator{\Prob}{\mathbf{Prob}}
\DeclareMathOperator{\WW}{\mathbf{W}}
\DeclareMathOperator{\AND}{AND}
\DeclareMathOperator{\OR}{OR}
\DeclareMathOperator{\Maj}{Maj}
\DeclareMathOperator{\Tribus}{Tribus}

\journame{La Gaceta de la RSME}
\yearofpublication{2025}
\volume{28}
\issuenumber{1}

\belongstopart{Artículos}

\title{Las funciones booleanas y el lema de Bonami}
\author{María José González, Paul MacManus y María Cristina Pereyra} 

\contact{María José González, Departamento de Matemáticas, Universidad de Cádiz}
{majose.gonzalez@uca.es}{}

\contact{Paul MacManus, Trevally Capital}
{paul.macmanus@gmail.com}{}

\contact{María Cristina Pereyra, Department of Mathematics and Statistics, University of New Mexico}
{crisp@math.unm.edu}{}

\begin{document}

\maketitle

\begin{abstract}
En este artículo expositivo estudiamos la relación entre las funciones booleanas y los teoremas de hipercontractividad de Aline Bonami. Nos concentramos en la teoría de la elección social, y presentamos algunos de los resultados más importantes en el área como los teoremas de Friedgut-Kalai-Naor (FKN) y de Kahn-Kalai-Linial (KKL), y la famosa conjetura \emph{entropía de Fourier / influencia}.
\end{abstract}


\section{Introducción}

La idea de escribir este artículo surge a partir de la conferencia impartida por Aline Bonami en el Seminari d’Anàlisi de Barcelona en la Universitat Autònoma de Barcelona el 29 de mayo de 2024, con ocasión de la celebración de la jubilación de Joaquim Bruna. El título de la conferencia era \emph{Riesz products, old and recent results}. Una vez finalizada su presentación, J. Bruna le preguntó por qué sus resultados sobre hipercontractividad son tan utilizados en el mundo aplicado. En este artículo expositivo intentaremos explicar esta
relación entre los teoremas de Bonami y su aplicación dentro de la teoría de la elección social.\footnote{\emph{Social choice theory} en la literatura matemática en inglés.}

\begin{figure}
   \centering
   \includegraphics[width=.8\textwidth]{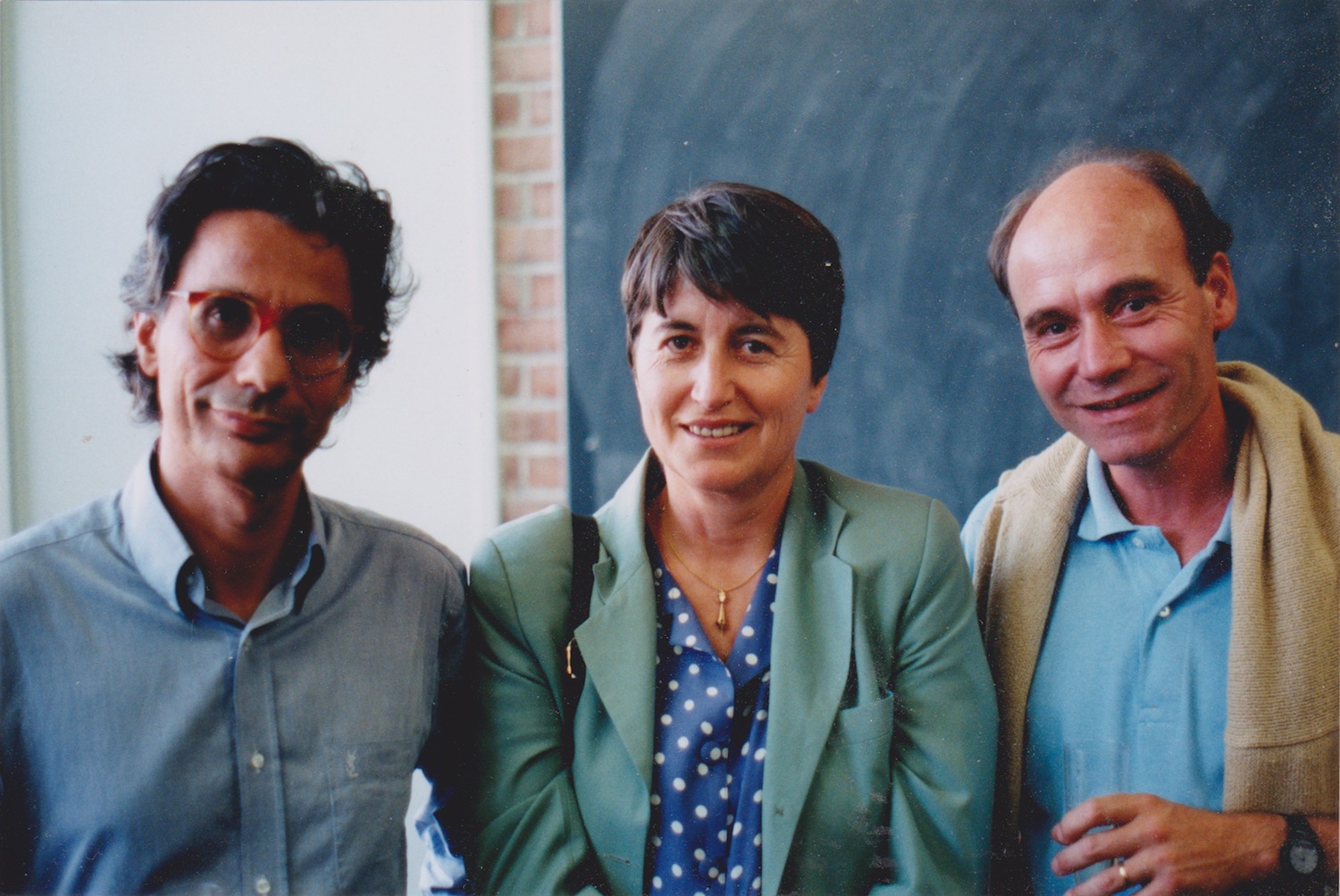}
   \pseudocaption{Aline Bonami con Nessim Sibony (izquierda) y Joaquim Bruna (derecha). Foto de Aline Bonami.}
\end{figure}

El puente de conexión son las funciones reales definidas en el cubo discreto $\{-1,1\}^n$, también llamado el \emph{cubo de Hamming}, prestando especial atención a las funciones booleanas: $f:\{-1,1\}^n\to\{-1,1\}$. Las funciones booleanas aparecen en todo tipo de contextos, tales como la teoría de computación, la combinatoria, la teoría de la elección social y la probabilidad discreta. También están íntimamente conectadas a las desigualdades funcionales y al análisis de Fourier en grupos discretos. El lema de Bonami al que se refiere el título de este artículo, y que será clave en nuestra exposición, forma parte de los resultados sobre hipercontractividad de las funciones reales en el cubo de Hamming demostrados por Bonami.

En los párrafos siguientes aparecerán varias nociones matemáticas con las que el lector quizás no esté familiarizado; esto no debe ser motivo de desánimo puesto que el objetivo de las secciones siguientes del artículo será presentar dichas nociones detalladamente.

En el contexto de la teoría de la elección social, una función booleana $f:\{-1,1\}^n\to\{-1,1\}$ representa un sistema de votación donde $n$ es el número de votantes. Cada votante elige entre dos candidatos A y B ($1$ y $-1$, respectivamente). La función $f$ define el sistema de votación, asignando $1$ o $-1$ a cada vector de unos y menos unos. Es decir, las entradas representan los votos individuales, y la salida representa el resultado de la votación. A modo ilustrativo consideremos un sistema electoral con dos votantes, en el que el candidato ganador es A a no ser que todos voten a B. En este caso la función booleana sería $f(1,1)=f(1,-1)=f(-1,1)=1$, $f(-1,-1)=-1$.

De forma más general consideramos funciones $f:\{-1,1\}^n\to\R$.
Asignando la medida uniforme en el cubo discreto $\{-1,1\}^n$, es decir, $\mu(x)=\frac{1}{2^n}$ para todo $x\in\{-1,1\}^n$, definimos la \emph{media} o \emph{esperanza} de $f$ como
\[
   \EE[f] =\frac{1}{2^n}\sum_{x\in\{-1,1\}^n} f(x),
\]
y la norma $\|f\|_p:= \big(\EE[|f|^p]\big)^{1/p}$ para todo $p\geq 1$. Cuando $p=2$, el espacio $L^2(\{-1,1\}^n)$ es un espacio de Hilbert, donde el producto interior se define como
\[
 \langle f,g\rangle := \EE[fg]= \frac{1}{2^n} \sum_{x\in \{-1,1\}^n}f(x)g(x).
\]
En este contexto se desarrolla la teoría del análisis de Fourier en el cubo discreto $\{-1,1\}^n$, que veremos detalladamente en la sección~\ref{sec:prelim}.

Denotaremos $[n]:=\{1,2,\dots, n\}$, notación común en el área de las funciones booleanas. Para cada subconjunto $S\subset [n]$, consideremos las funciones
\begin{align*}
   \chi_{\emptyset} &= 1,
   \\
   \chi_S(x_1,\dots,x_n) &= \prod_{i\in S} x_i.
\end{align*}
Como veremos más adelante, las funciones $\chi_S$, $S\subset [n]$, forman un sistema ortogonal, y toda función $f:\{-1,1\}^n\to\R$ se puede escribir de forma única como un polinomio multilineal:\footnote{La representación como polinomio multilineal es única pero no la representación como polinomio, puesto que, por ejemplo, $x_i^3 = x_i$.}
\begin{equation}\label{ecn:desarrollo_fourier}
   f(x_1,x_2,\dots, x_n) = \sum_{S\subset [n]} a_S \chi_S,
\end{equation}
con $a_S\in \R$.
Por ejemplo, la función $f$ definida anteriormente tiene representación de la forma:
\[
f(x_1,x_2)=1-\frac{1}{2}(1-x_1)(1-x_2)=\frac{1}{2}(1+x_1+x_2-x_1x_2).
\]

Las funciones $\chi_S$ representan la \emph{base de Fourier}, y $a_S=\langle f,\chi_S\rangle$ los \emph{coeficientes de Fourier} de la función $f$, que a partir de ahora denotaremos $\widehat{f}(S)=a_S$. De esta forma, toda función real definida en el cubo de Hamming se puede expresar en su desarrollo de Fourier \eqref{ecn:desarrollo_fourier} y el grado de $f$ se define como el grado de este polinomio.

El \emph{lema de Bonami} establece que cuando $f:\{-1,1\}^n\to\R$ tiene grado $d$, se verifica la siguiente desigualdad de tipo Hölder inversa:
\[
\|f\|_4\leq 3^{d/2}\|f\|_2,
\]
es decir, $\EE[f^4]\leq 9^d\big(\EE[f^2]\big)^2$.

Para explicar la relación de este resultado con las funciones booleanas necesitamos introducir uno de los conceptos más importantes en el análisis de las funciones booleanas: la \emph{influencia}.
Regresando a la interpretación de una función booleana
$f:\{-1,1\}^n\to\{-1,1\}$ como un sistema de votación, resulta natural definir la influencia del votante $i$ como el porcentaje de escenarios en los que el voto de~$i$ es decisivo, es decir, si este cambia su voto de $1$ a $-1$ o viceversa, el resultado de la votación varía también. En el ejemplo anterior, el voto del votante $1$ es decisivo en los escenarios $(1, -1)$ y $(-1, -1)$. Por tanto su influencia es $1/2$, igual que la del votante $2$.

En general, dada una función booleana $f$, la \emph{influencia de la variable} $i$ se define como la probabilidad de que $f(x)\neq f(x^i)$, donde $x^i=(x_1,\dots, -x_i,\dots, x_n)$, es~decir,
\[
\II_i(f) = \frac{1}{2^n}\left|\{x\in\{-1,1\}^n : f(x)\neq f(x^i)\}\right|.
\]
La \emph{influencia total} de la función booleana $f$ se define como $\II(f)=\sum_{i=1}^n \II_i(f)$. En términos de los coeficientes de Fourier
\[
   \II_i(f)= \sum_{S\subset [n] \colon i\in S} |\widehat{f}(S)|^2,
\]
e
\[
   \II(f)=\sum_{S\subset [n]} |\widehat{f}(S)|^2 |S|
\]
(ver la sección~\ref{sec:influence}), donde $|S|$ es el tamaño del conjunto $S$.

Puesto que $\EE(\chi_S)=0$ para todo $S\neq\emptyset$, la esperanza y varianza de $f$ vienen dadas por $\EE[f]= \widehat{f}(\emptyset)$ y
\[
   \Var(f)=\sum_{S\subset [n] \colon S\neq \emptyset} |\widehat{f}(S)|^2
\]
respectivamente. Observemos ahora que $\II(f)\geq \Var(f)$. Esta desigualdad nos permite obtener de forma trivial la siguiente cota superior para la máxima influencia de~$f$: $\max_{i\in [n]} \II_i(f) \geq \frac{1}{n} \Var(f)$. Sin embargo, el sorprendente \emph{teorema de Kahn-Kalai-Linial} (KKL) establece que
\[
\max_{i\in [n]} \II_i(f) \geq C \frac{\log n}{n} \Var(f),
\]
donde $C>0$ es una constante absoluta. Este es uno de los resultados más relevantes en el análisis de las funciones booleanas. Como veremos en la sección~\ref{sec:social-choice-meets-Bonami}, el ingrediente clave en su demostración es el lema de Bonami. Una consecuencia del teorema de KKL es que en el caso de un sistema de votación monótono con $n$ votantes, existe una coalición de tamaño comparable a $n/\log n$ que básicamente controla el resultado de la votación. Este resultado es conocido como \emph{el problema de Ben-Or y Linial}.

El lema de Bonami es equivalente a un caso particular del llamado \emph{teorema de hipercontractividad}, también un resultado de Bonami \cite[capítulo III, teorema~3]{Bon70}, y al que dedicaremos la sección~\ref{sec:hypercontractitvity}. Este teorema dice que ciertos multiplicadores de Fourier $T_\rho$ sobre funciones $f:\{-1,1\}^n\to\R$ son contracciones y regularizaciones de~$f$, en el sentido de que $\|T_{\rho} f\|_q\leq \|f\|_p$ para un rango de $q\geq p$. Estas desigualdades tienen aplicaciones importantes en los campos de la ciencia de la computación y la teoría de grafos aleatorios \cite{KaS05}.

Finalizamos esta introducción enunciando una de las conjeturas más famosas en el área, que analizaremos con más detalle en la sección~\ref{sec:EntropyInfluence}. La identidad de Parseval garantiza que para funciones booleanas $f:\{-1,1\}^n\to \{-1,1\}$ tenemos que $\sum_{S\subset [n]} |\widehat{f}(S)|^2=1$ (ver la sección \ref{sec:prelim}). Por tanto, podemos definir una medida de probabilidad en $\{S:S\subset [n]\}$ que asocia a cada subconjunto $S$ la probabilidad $|\widehat{f}(S)|^2$. Denotando la entropía de esta distribución de probabilidad por $\HH(f)$, el problema que permanece abierto es la \emph{conjetura entropía de Fourier / influencia}, que abreviaremos \emph{conjetura EFI}: existe $C>0$ tal que, para toda función booleana~$f$, $\HH(f)\leq C \II(f)$.

Es decir, para funciones booleanas $f:\{-1,1\}^n\to\{-1,1\}$,
\[
\sum_{S\subset [n]} |\widehat{f}(S)|^2\log_2 \frac{1}{|\widehat{f}(S)|^2}
\leq C \sum_{S\subset [n]} |\widehat{f}(S)|^2|S|.
\]
La importancia de esta conjetura radica en el hecho de que, de ser cierta, implicaría resultados como el teorema de KKL o la llamada conjetura de Mansur \cite{Kal07}.

Este artículo está organizado de la siguiente manera:

En la sección~\ref{sec:prelim} presentamos las definiciones básicas de la teoría de funciones booleanas, la descomposición de Fourier y las nociones de influencia individual y total. Discutimos los ejemplos canónicos en el contexto de la teoría de la elección social, como la función mayoría, las dictaduras, las funciones $\AND$ y $\OR$, las juntas, las tribus, y las funciones booleanas monótonas.

En la sección~\ref{sec:Bonami} demostramos el lema de Bonami y deducimos algunos corolarios que nos serán útiles más adelante. Enunciamos el teorema de hipercontractividad y algunos casos particulares. Entre ellos, el llamado teorema de hipercontractividad-$(4,2)$, que es equivalente al lema de Bonami.

En la sección~\ref{sec:social-choice-meets-Bonami} demostramos que una función booleana afín es una constante o una dictadura. Como consecuencia del lema de Bonami mostramos los teoremas de Friedgut-Kalai-Naor (FKN) y de Kahn-Kalai-Linial (KKL). El teorema de FKN muestra que una función que es aproximadamente lineal debe estar cerca de una dictadura. Para terminar la sección, presentamos una aplicación del teorema de KKL a la resolución del problema de Ben-Or y Linial.

En la sección~\ref{sec:EntropyInfluence} discutimos las conjeturas EFI y MEFI (\emph{mínima entropía de Fourier / influencia}). Mostramos cómo obtener una primera aproximación a la conjetura EFI pero con una cota dependiente de la dimensión $n$ y de orden $\log n$. También estudiamos las conjeturas en determinadas clases de funciones, donde o bien se han comprobado o se han obtenido resultados parciales. Finalmente mostramos la conexión entre el teorema de KKL y las conjeturas MEFI y EFI.

Los autores desean expresar su agradecimiento a Joaquim Bruna, de la Universitat Autònoma de Barcelona, por sus comentarios y correcciones sobre una versión previa de este artículo, y a José Luis Fernández, de la Universidad Autónoma de Madrid, por las múltiples y valiosas conversaciones sobre la conjetura EFI. Agradecemos también al evaluador/a por la exhaustiva lectura de este artículo y por sus comentarios, que han ayudado a mejorar notablemente esta presentación.

\section{Funciones booleanas, elección social, e influencia}\label{sec:prelim}

En esta sección, presentamos las funciones booleanas, su desarrollo de Fourier y el concepto de la influencia. Miramos también algunos ejemplos canónicos de la teoría de la elección social, en concreto, ejemplos de sistemas de votación. Una referencia excelente en esta materia es el libro \emph{Analysis of Boolean functions} por Ryan O'Donnell \cite{O'Do21}; en particular, las notas históricas al final de cada capítulo son muy informativas.

\subsection{Cubo de Hamming, funciones booleanas, y caracteres}

El \emph{cubo de Hamming de dimensión $n$}, escrito como $\{-1,1\}^n$, consiste en vectores de longitud $n$ con valores en $\pm 1$. Es decir:
\[
\{-1,1\}^{n}=\{ (x_1,x_2,\dots,x_n): x_i=\pm 1\}.
\]
Observemos que este conjunto tiene $2^n$ elementos.

El conjunto de funciones definidas en el cubo de Hamming de dimensión $n$ con valores reales, $f:\{-1,1\}^{n} \to \R$, es un espacio vectorial que denotamos $\mathcal{F}_n$. Una función $f:\{-1,1\}^{n} \to \{-1,1\}$ se llama \emph{función booleana} y el conjunto de estas funciones se escribe $\mathcal{B}_n$. En total, hay $2^{2^n}$ funciones booleanas definidas en $\{-1,1\}^n$.

Asignamos la medida uniforme al cubo de Hamming. Bajo esta medida, las variables $x_i$ son variables aleatorias independientes uniformemente distribuidas que toman los valores $1$ y $-1$. Los elementos $f(x_1,\dots, x_n)$ de $\mathcal{F}_n$ representan por tanto variables aleatorias.

Definimos el producto interior de dos funciones $f, g \in \mathcal{F}_n$ por:
\[
\langle f,g\rangle := \EE[fg]= \frac{1}{2^n} \sum_{x\in \{-1,1\}^{n}}f(x)g(x).
\]
La norma correspondiente es $\|f\|_2=\sqrt{\langle f,f\rangle}$ y el espacio de Hilbert resultante $(\mathcal{F}_n, \langle \cdot , \cdot \rangle)$ se escribe $L^2(\{-1,1\}^n)$. Todas las funciones booleanas tienen norma igual a~$1$.

Para $S\subset [n]:=\{1,2,\dots ,n\}$, el \emph{carácter}, o \emph{función de paridad}, correspondiente se define como:
\begin{align*}
 \chi_{\emptyset}(x_1,\dots,x_n) &:=1,
 \\
 \quad \chi_S(x_1,\dots,x_n) &:=\prod_{i\in S} x_i \text{ para todo } S\neq \emptyset.
\end{align*}
El número de caracteres es igual al número de subconjuntos $S\subset [n]$, es decir, $2^n$. Vamos a ver a continuación que la colección de caracteres es una base ortonormal para $L^2(\{-1,1\}^n)$. Nótese que cada carácter es una función booleana.

Los caracteres cumplen:
\begin{align*}
\chi_S(x)\chi_U(x)&=\chi_{S\Delta U}(x),
\\
\EE[\chi_V] &= \delta_{\emptyset}(V),
\end{align*}
donde $S\Delta U$ es la diferencia simétrica de $S$ y $U$, y donde $\delta_{\emptyset}(V)=1$ si $V=\emptyset$ y $0$ en caso contrario. Tomando $V=S\Delta U$ se obtiene:
\[
\langle \chi_S,\chi_U\rangle = \EE[\chi_S \chi_U] =
\EE[\chi_{S\Delta U}] = \delta_{\emptyset}(S\Delta U).
\]
Como $S \Delta U=\emptyset$ si, y solo si, $S=U$, hemos demostrado que los caracteres forman un conjunto ortonormal en $L^2(\{-1,1\}^n)$. ¿Y por qué forman una base? Porque $L^2(\{-1,1\}^n)$, visto como espacio vectorial, tiene dimensión $2^n$. Esto último se debe a que las funciones indicatrices, $\mathbbm{1}_{\{a\}}(x)=\delta_a(x)$ para $a\in \{-1,1\}^n$, forman una base ortogonal para $L^2(\{-1,1\}^n)$.

Alternativamente, para la función $\mathbbm{1}_{\{a\}}$, $a\in \{-1,1\}^n$, tenemos:
\begin{align}
   \mathbbm{1}_{\{a\}}(x)
   &= \frac{1+a_1x_1}{2} \cdot
   \frac{1+a_2x_2}{2} \cdot
   \ldots \cdot \frac{1+a_nx_n}{2}
   \notag
   \\
   &= \frac{1}{2^n} \sum_{S\subset [n]} \Big( \prod_{i\in S}a_i \Big) \chi_S(x). \label{eq:expansion-standard-basis-in-parity-functions}
\end{align}
De modo que cualquier función real es una combinación lineal de caracteres y, como estos son ortonormales, tienen que ser una base ortonormal para $L^2(\{-1,1\}^n)$.

\subsection{Coeficientes de Fourier}\label{sec:coeficientes-Fourier}

Como los caracteres forman una base ortonormal, cada función $f:\{-1,1\}^{n}\to \R$ tiene un desarrollo único del estilo:
\begin{equation}
\label{desarrollo_fourier}
f(x) = \sum_{S\subset [n]} \widehat{f}(S) \chi_S(x),
\end{equation}
donde $\widehat{f}(S)=\langle f,\chi_S\rangle$ representa el \emph{coeficiente de Fourier de $f$} asociado a $S$ y el polinomio~\eqref{desarrollo_fourier} se denomina el \emph{desarrollo de Fourier de $f$}.

Por ejemplo, por \eqref{eq:expansion-standard-basis-in-parity-functions} los coeficientes de Fourier de $\mathbbm{1}_{\{a\}}$ son:
\[
   \widehat{\mathbbm{1}_{\{a\}}}(S)
   = \frac{1}{2^n} \prod_{i\in S}a_i
   = \frac{1}{2^n} \chi_S(a).
\]

Como consecuencia del desarrollo de Fourier \eqref{desarrollo_fourier}, tenemos la \emph{identidad de Plancherel}
\[
   \EE[fg] = \frac{1}{2^n}\sum_{x\in \{-1,1\}^n} f(x) g(x)
   = \sum_{S\subset [n]} \widehat{f}(S) \widehat{g}(S)
\]
y la \emph{identidad de Parseval}
\begin{equation}
\label{Parseval}
   \EE[f^2] = \frac{1}{2^n} \sum_{x\in \{-1,1\}^n} |f(x)|^2
   = \sum_{S\subset [n]} |\widehat{f}(S)|^2.
\end{equation}
En particular, si $f$ es booleana, tenemos que $\sum_{S\subset [n]} |\widehat{f}(S)|^2 = 1$.

Además, para toda $f:\{-1,1\}^{n} \to \R$, se pueden escribir cantidades como la \emph{esperanza}, la \emph{varianza} y la \emph{covarianza} en términos de los coeficientes de Fourier:
\begin{align*}
   \EE[f] &=\widehat{f}(\emptyset),
   \\
   \Var[f] &=\EE\big[f^2\big] -\EE[f]^2
   =\sum_{S\subset [n], S\neq\emptyset} |\widehat{f}(S)|^2,
   \\
   \Cov[f,g] &= \EE\big[fg \big]-\EE[f] \EE[g]
   = \sum_{S\subset [n], S\neq\emptyset} \widehat{f}(S) \widehat{g}(S).
\end{align*}

\subsection{Estructura de grupos}

El cubo de Hamming es un grupo multiplicativo de $2^n$ elementos, con el producto de dos elementos definido por $x\cdot y = (x_1y_1,x_2y_2,\dots,x_ny_n)$. Evidentemente, la identidad es el punto $\mathbf{1}:=(1,1,\dots, 1)$ y cada elemento es su propio inverso. Con esta estructura el grupo es isomorfo al grupo aditivo $\Z_2^n$.

También forman un grupo multiplicativo las funciones booleanas, con el producto de dos funciones $f$ y $g$ definido por $(f\cdot g)(x) = f(x)g(x)$. Con esta estructura, $\mathcal{B}_n$ es isomorfo al grupo aditivo $\Z_2^{2^n}$.

El grupo dual de $\{-1,1\}^{n}$ consiste en las funciones $\phi: \{-1,1\}^{n} \to \mathbb{T}$ (los números complejos de norma $1$) que respetan la estructura de grupo: $\phi(xy)=\phi(x)\phi(y)$. Es decir, el grupo dual consiste en los homomorfismos que van del cubo de Hamming a $\mathbb{T}$. En particular, $\phi(\mathbf{1})=1$ y $\phi(x)^2=\phi(x^2)=\phi(\mathbf{1})=1$ y, por tanto, todos los elementos del grupo dual son funciones booleanas.

Para $S\subset [n]$, $\chi_S$ es booleana y
\[
   \chi_S(xy)
   =\prod_{i\in S} (xy)_i
   = \prod_{i\in S} (x_iy_i)
   = \prod_S x_i  \prod_S y_i
   = \chi_S(x) \chi_S(y),
\]
de modo que $\chi_S$ es un elemento del grupo dual del cubo de Hamming. ¡Resulta que los únicos elementos del grupo dual son los caracteres! Se puede consultar por ejemplo \cite[sección 8.5]{O'Do21}, o cualquier texto de análisis de Fourier en grupos, como el texto clásico de Rudin~\cite[capítulo 1]{Ru62}.

\subsection{Ejemplos canónicos}\label{sec:canonical-examples}

Aquí presentamos algunos ejemplos relacionados con la teoría de la elección social, donde se investigan sistemas de votación o de agregación de opiniones. Digamos que tenemos dos alternativas $1$ y $-1$ y una sociedad de $n$ votantes. Cada punto del cubo de Hamming representa un conjunto de preferencias individuales de todos los miembros de la sociedad. Una función booleana es entonces una regla que decide una preferencia del conjunto de la sociedad a partir de las preferencias individuales.

\subsubsection{Dictadura, mayoría, \texorpdfstring{$\AND$}{AND} y \texorpdfstring{$\OR$}{OR}}

Las funciones dictadura son los caracteres que corresponden a los conjuntos de un solo elemento, los $\chi_{\{i\}}$.

\begin{example}
La i-ésima función dictadura $d_i:\{-1,1\}^n\to\{-1,1\}$ se define como $d_i(x):=x_i$.
\end{example}

En este caso, el resultado de la elección depende únicamente del votante $i$, el dictador. Es evidente que los coeficientes de Fourier de $d_i$ vienen dados por:
\[
   \widehat{d_i}(S) = \begin{cases}
   1 &\text{si } S = \{i\},
   \\
   0 &\text{en caso contrario}.
   \end{cases}
\]

Entre los mecanismos de votación más comunes se encuentra el de la mayoría.

\begin{example}
Para $n$ impar,\footnote{Cuando $n$ es par, se puede usar la misma fórmula cuando $x_1+\dots +x_n\neq 0$ para definir la función $\Maj_n(x)$ y cuando la suma es 0 se puede asignar valor $1$ o $-1$ arbitrariamente. En la práctica, cuando hay un empate, una autoridad superior decide, por ejemplo, un tribunal supremo.} la función mayoría $\Maj_n: \{-1,1\}^n\to \{-1,1\}$ se define como $\Maj_n(x)=\sgn (x_1+x_2+\dots +x_n)$.
\end{example}

Esta función toma el valor $1$ si, y solo si, el número de coordenadas que toman el valor $1$ es mayor que el número de las que toman el valor $-1$.

Calcular los coeficientes de Fourier de esta función es algo complejo, pero tenemos la fórmula siguiente \cite[teorema 5.19]{O'Do21}:
\[
   \widehat{\Maj}_n(S)
   = \begin{cases}
   0 &\text{si } |S| \text{ es par},
   \\
   (-1)^k \frac{\binom{(n-1)/2}{k}}{\binom{n-1}{2k}}
   \frac{2}{2^n} \binom{n-1}{(n-1)/2}
   &\text{si } |S|=2k+1.
   \end{cases}
\]

\begin{example}
La función $\OR_n:\{-1,1\}^n\to\{-1,1\}$ se define como
\[
   \OR_n (x) = \begin{cases}
   +1 &\text{si } x \neq (-1,-1,\dots, -1),
   \\
   -1 &\text{si } x = (-1,-1,\dots, -1).
   \end{cases}
\]
\end{example}

Esta función corresponde al caso en el que el candidato $1$ gana salvo que todos los votantes voten en su contra. Se puede escribir también como $\OR_n(x) = \max(x_i)$ o como $\OR_n (x) = 1-2\mathbbm{1}_{\{\mathbf{-1}\}}(x)$. Como ya sabemos los coeficientes de Fourier de la función indicatriz $\mathbbm{1}_{\{\mathbf{-1}\}}$, obtenemos:
\[
   \widehat{\OR_n}(S) = \begin{cases}
   1-1/2^{n-1} &\text{si } S=\emptyset,
   \\
   (-1)^{|S|+1}/2^{n-1} &\text{en caso contrario}.
   \end{cases}
\]

\begin{example}
La función $\AND_n:\{-1,1\}^n\to\{-1,1\}$ se define como
\[
   \AND_n(x) = \begin{cases}
   -1 &\text{si } x \neq (1,1,\dots, 1),
   \\
   +1 &\text{si } x = (1,1,\dots, 1).
   \end{cases}
\]
\end{example}

En este caso el candidato $1$ gana solo cuando todos votan a su favor. Esta función se puede escribir como $\AND_n(x)=\min(x_i)$. También, tenemos $\AND_n(x) = -\OR_n(-x)$; por tanto los coeficientes de Fourier de esta función son:
\[
   \widehat{\AND_n}(S) = \begin{cases}
   -1+1/2^{n-1} &\text{si } S=\emptyset,
   \\
   1/2^{n-1} &\text{en caso contrario}.
   \end{cases}
\]

\subsubsection{Juntas, tribus, y funciones monótonas}

Los siguientes son ejemplos algo más elaborados que los anteriores.

\begin{example}
Una función booleana es una $k$-junta si depende solo de $k$ de sus $n$ coordenadas. Es decir, existen $x_{i_1},x_{i_2},\dots, x_{i_k}$ y una función booleana $g : \{-1,1\}^k \to \{-1,1\}$ tal que $f(x) =g(x_{i_1},x_{i_2},\dots, x_{i_k})$.
\end{example}

La función dictadura corresponde al caso $k=1$ de la función junta.
Los coeficientes de Fourier de una $k$-junta $f$ vienen dados por $\widehat{f}(S)=\widehat{g}(S)$ si $S\subset \{i_1,i_2,\dots, i_k\}$ y $\widehat{f}(S)=0$ en caso contrario.

En el siguiente ejemplo, los votantes se dividen en \emph{tribus} disjuntas y el candidato~$1$ gana si, y solo si, por lo menos una de las tribus vota unánimemente a su favor.

\begin{example}\label{eg:tribes}
Para una partición $\mathcal{P}=\{I_1, I_2, \dots, I_s\}$ de $[n]$, la función tribus de $\mathcal{P}$ se define como
\[
\Tribus_{\mathcal{P}}(x)= \OR_s\big(\AND(x^{(1)}),\dots \AND(x^{(s)})\big),
\]
donde $x^{(k)}$ son los elementos de $x$ con índices en $I_k$.
\end{example}

Las funciones tribus son importantes como ejemplos para mostrar que resultados fundamentales, como el teorema de Kahn–Kalai–Linial (KKL), son óptimos. Cuando la partición $\mathcal{P}$ es uniforme (los $I_k$ tienen el mismo número de elementos), se pueden calcular los coeficientes de Fourier de la función tribus. Véase~\cite[proposición~4.14]{O'Do21}.

Las funciones booleanas descritas hasta el momento, salvo las juntas, tienen una propiedad importante en común: son monótonas.

\begin{definition}
Una función booleana $f:\{-1,1\}^n\to \{-1,1\}$ es monótona si $f(x) \leq f(y)$ cuando $x_i\leq y_i$ para todo $i\in [n]$.
\end{definition}

Nótese que una junta es monótona si, y solo si, la función $g$ que la define es monótona.

En términos de la teoría de la elección social, una función monótona garantiza que un cambio de voto a favor de un candidato no perjudica a dicho candidato en el resultado de la elección.

\subsection{La influencia}\label{sec:influence}

Si $f:\{-1,1\}^n\to\{-1,1\}$ representa un sistema de votación, la influencia del votante $i$ es una medida del impacto electoral de un cambio de voto por parte del votante $i$. Con este fin, definimos $x^{(i)}$ como el vector que es igual a $x$ salvo en la $i$-ésima coordenada, donde tiene un cambio de signo.

\begin{definition} La influencia de la coordenada $i$ en una función booleana $f:\{-1,1\}^n \linebreak\to\{-1,1\}$ es
\[
\II_i(f)= \Prob\big\{f(x)\neq f(x^{(i)})\big\}.
\]
\end{definition}

Así que la influencia es la probabilidad de que un cambio de voto por parte del votante~$i$ cambie el resultado de la elección. Un punto $x$ en el que $f(x)\neq f(x^{(i)})$ se llama un \emph{pivote} de $f$ para la coordenada $i$.

\begin{definition} La influencia total de una función booleana $f:\{-1,1\}^n\to \{-1,1\}$~es
\[
\II(f)= \sum_{i\in [n] } \II_i(f).
\]
\end{definition}

Por tanto, la influencia total es una medida de la sensibilidad de $f$ a cambios de voto.
Veremos pronto cómo extender estas definiciones al caso más general de funciones reales definidas en el cubo de Hamming.

La influencia tiene un vínculo estrecho con la derivada de la función. En general, para toda función $g:\{-1,1\}^n\to \R$ la \emph{derivada $i$-ésima} se define como
\[
   D_ig({x}) = \frac{g({x}^{(i\to1)}) - g({ x}^{(i\to-1)})}{2},
\]
donde usamos la notación ${x}^{(i\to b)}=(x_1,\dots,x_{i-1},b,x_{i+1},\dots, x_n)$.

Las derivadas actúan como derivadas formales en los caracteres, es decir,
\[
   D_i\chi_S = \begin{cases}
   \chi_{S\setminus\{i\}} &\text{si } i\in S,
   \\
   0 &\text{en caso contrario.}
   \end{cases}
\]
Por tanto, por linealidad, también son derivadas formales de los desarrollos de Fourier. Formalmente, y como cada función $g:\{-1,1\}^n\to \R$ tiene su desarrollo de Fourier \eqref{desarrollo_fourier}, tenemos
\begin{equation}\label{eqn:derivativeFourier}
   D_ig(x)=\sum_{S\subset [n] \colon S\ni i}\widehat{g}(S) \chi_{S\setminus \{i\}}(x).
\end{equation}

En el caso de una función booleana, las derivadas solo toman valores en $\{-1,0,1\}$, y toman el valor $0$ solo en los puntos que no son pivotes de $f$. Por tanto, para una función booleana $f$, se verifica
\begin{equation}
\label{eqn:influencia-derivada}
   \II_i(f) = \EE[D_if(x)^2] = \|D_i f\|_2^2
\end{equation}
e
\[
   \II(f) = \sum_{i\in [n]} \II_i(f)
   = \sum_{i\in [n]} \|D_i f\|_2^2 =: \EE[\nabla f]^2.
\]
Aquí, $\nabla f=(D_1f,\dots,D_nf)$ es el gradiente de $f$.

Por la identidad de Parseval~\eqref{Parseval}, concluimos que
\begin{equation}\label{eq:influencia_i}
\II_i(f)= \sum_{S\subset [n] \colon S\ni i}|\widehat{f}(S)|^2.
\end{equation}

\begin{proposition}[influencia y Fourier]
\label{prop:influence-Fourier}
Para una función booleana
$f:\{-1,1\}^n\to \{-1,1\}$, tenemos
\[
\II(f)=\sum_{S\subset [n]} |S| |\widehat{f}(S)|^2,
\]
donde recordamos que $|S|$ denota el tamaño del conjunto $S$.
\end{proposition}

\begin{proof}
La influencia total es la suma de las influencias individuales y para cada $S$ el término $|\widehat{f}(S)|^2$ aparece en exactamente $|S|$ de las influencias individuales (véase la fórmula (\ref{eq:influencia_i})).
\end{proof}

Observemos que esta proposición y la identidad \eqref{eq:influencia_i} nos permiten definir la influencia total y la influencia individual de cualquier función $g:\{-1,1\}^n \to \R$ en términos de sus coeficientes de Fourier:
\[
\II(g)=\sum_{S\subset [n]} |S| |\widehat{g}(S)|^2
\quad\text{e} \quad
\II_i(g)= \sum_{S\subset [n]:S\ni i}|\widehat{g}(S)|^2.
\]

Recordemos que, por la identidad de Parseval \eqref{Parseval}, para toda función booleana se cumple $\sum_{S\subset [n}|\widehat{f}(S)|^2=1$. Por tanto, podemos interpretar que los coeficientes de Fourier $|\widehat{f}(S)|^2$ representan una medida de probabilidad sobre los subconjuntos $S\subset [n]$. Con este enfoque, la proposición anterior dice que la influencia total es la esperanza de $|S|$ bajo esta medida.

En el caso de las funciones booleanas monótonas, la relación entre la influencia y los coeficientes de Fourier se simplifica: véase \cite[proposición 2.21]{O'Do21}.

\begin{proposition}\label{prop:influenceMonotone}
Si una función booleana $f:\{-1,1\}^n\to \{-1,1\}$ es monótona, entonces $\II_i(f)=\widehat{f}(\{i\})$.
\end{proposition}

\begin{proof}
Como $x^{(i\to -1)}\leq x^{(i\to1)}$ por monotonicidad, $D_if(x)$ solo toma valores en $\{0,1\}$. Además, $D_if(x)=1$ si, y solo si, $x$ es un pivote de $f$ para la coordenada $i$. Así que,
\[
\II_i(f) = \EE [ D_if]= \widehat{D_if}(\emptyset)=\widehat{f}(\{i\}).
\]
La última igualdad es consecuencia de la relación \eqref{eqn:derivativeFourier}.
\end{proof}

\begin{corollary}
La influencia total de una función monótona $f:\{-1,1\}^n\to \{-1,1\}$ viene dada por
\[
\II(f) = \sum_{i=1}^n \widehat{f}(\{i\}).
\]
\end{corollary}

Terminamos esta sección con un análogo para funciones definidas en el cubo de Hamming de la desigualdad de Poincaré en el entorno continuo. La desigualdad de Poincaré permite controlar la variación de una función por la norma de su derivada. Por ejemplo, en un intervalo acotado en la recta, $I$, y para $f,f' \in L^2(I)$, tenemos
\[
   \int_I |f(x)-f_I|^2 \, dx
   = \int_I \big(f^2(x)-f_I^2\big) \, dx
   \leq \frac{|I|^2}{\pi^2} \int_I |f'(x)|^2 \, dx,
\]
donde $f_I=\frac{1}{|I|} \int_I f(x)\,dx$. Esta variante de la desigualdad se debe a Wirtinger. Nótese que la función $f(x)=\cos \big(\pi (x-a_I)/|I|\big)$, donde $I=[a_I,b_I]$, demuestra que la desigualdad es óptima. Para informarse sobre esta y muchas otras desigualdades puede consultarse el clásico \cite[sección 7.7]{HLP23}.

\begin{proposition}[desigualdad de Poincaré]
\label{thm:Poincare's Inequality}
Para toda función $f:\{-1,1\}^n \to \R$, se cumple
\[
   \Var(f) \leq \II(f)
\]
o, equivalentemente,
\[
   \EE[f^2-\EE(f)^2]\leq \EE[\nabla f]^2.
\]
\end{proposition}

\begin{proof} Tenemos
\[
   \Var(f) = \sum_{S\subset [n],\, S\neq \emptyset} |\widehat{f}(S)|^2
   \leq \sum_{S\subset [n]} |S|\, |\widehat{f}(S)|^2 = \II(f).
\]
La primera igualdad es la fórmula de Fourier para la varianza (ver sección \ref{sec:coeficientes-Fourier}) y la segunda es la definición de la influencia para funciones definidas en el cubo de Hamming con valores reales.
\end{proof}

Nótese que la función $f(x)=x_1$ demuestra que la desigualdad es óptima para cualquier $n$.

\section{El lema de Bonami y la hipercontractividad}\label{sec:Bonami}

En esta sección presentamos el lema de Bonami y, de forma más general, el teorema de hipercontractividad, también un resultado de Bonami de finales de la década de los~60~\cite{Bon68, Bon70}. Es importante resaltar también los resultados unos años más tarde de Beckner~\cite{Bec75}, que obtiene, independientemente, desigualdades óptimas de hipercontractividad, en particular la constante óptima en la desigualdad de Hausdorff-Young. Además de Bonami y Beckner, la historia de los teoremas de hipercontractividad involucra a muchos otros distinguidos analistas y probabilistas, entre ellos Gross, Paley y Talagrand \cite{Gro75,Pal32,Tal93}. Recomendamos las notas históricas en \cite[pp. 278--281]{O'Do21} para más detalles y referencias.

El blog de Kalai \cite{Kal20a} atribuye a Simon haber acuñado la palabra hipercontractividad en los años 70. Esto lo confirma el propio Simon en~\cite{DGS90}, una nota histórica acompañada de una larga bibliografía, donde se cuenta cómo las ideas de hipercontractividad han tenido mucha influencia en análisis con ramificaciones a la teoría cuántica de campos, a los operadores de Schrödinger y a la teoría de semigrupos hipercontractivos.
Nosotros nos centraremos en las conexiones con las funciones booleanas.

\subsection{El Lema de Bonami}

Recordemos que, para $1\leq p < \infty$, la norma en $L^p$ de una función $f:\{-1,1\}^n\to \R$ se define como
\[
\|f\|_p = \left( \EE[|f|^p] \right)^{1/p},
\]
y que el grado de $f$ es el grado del polinomio dado por su desarrollo de Fourier~\eqref{desarrollo_fourier}.

En esta sección discutimos el celebrado lema de Bonami sobre funciones reales definidas en el cubo de Hamming.
Este resultado es interpretado por los probabilistas en términos de lo que denominan variables aleatorias \emph{razonables}. Dado $B\geq 1$, diremos que una variable aleatoria $X$
es $B$-razonable si $E(X^4)\leq B (E(X^2))^2$, o equivalentemente, si $\|X\|_4\leq B^{1/4} \|X\|_2$. Nótese que $\|X\|_2 \leq \|X\|_4$ no es más que la desigualdad de Hölder para una medida de probabilidad.

En estadística, el parámetro $B$ está relacionado con la \emph{curtosis} de $X$, que indica el grado de concentración de los valores de $X$ alrededor de la media. Las colas de las variables $B$-razonables son finas, en el sentido de que satisfacen la desigualdad
\[
\Prob\big\{ |X|\geq t \|X\|_2 \big\}\leq \frac{B}{t^4}, \quad t>0,
\]
a la vez que los valores de $X$ no están demasiado concentrados alrededor de $0$,
\[
\Prob \big\{ |X|\geq t \|X\|_2 \big\} \geq \frac{(1-t^2)^2}{B}, \quad 0\leq t \leq 1.
\]
La primera desigualdad es consecuencia de la desigualdad de Markov aplicada a la función $|X|^4$, y la segunda se deriva de la desigualdad de Paley-Zygmund para $|X|^2$. Véase
\cite[sección 9.1]{O'Do21}, por ejemplo.

El lema de Bonami nos dice que los polinomios en el cubo de Hamming son razonables.

\begin{lemma}[lema de Bonami]\label{lem:Bonami's}
Sea $f:\{-1,1\}^n\to \R$ un polinomio de grado $d\geq 0$. Entonces
\[
\|f\|_4\leq 3^{d/2}\|f\|_2,
\]
es decir,
\[
\EE[f^4]\leq 9^d \EE[f^2]^2.
\]
\end{lemma}

Observemos que, como $x_i^{2k}=1$, $x_i^{2k+1}=x_i$, cualquier polinomio se puede reducir a una combinación lineal de caracteres y, por tanto, se puede reducir a un polinomio de grado no mayor que $n$.

\begin{proof}[Demostración del lema de Bonami.] Procederemos por inducción en el número de variables $n$ (la dimensión del cubo de Hamming). Cuando $n=1$ entonces $f(x_1)=a_0+a_1x_1$ (puesto que no hay otros caracteres). Calculemos el cuadrado y la cuarta potencia de $f$:
\begin{align*}
   f^4(x_1) &= a_0^4+a_1^4+6a_0^2a_1^2+4a_0a_1(a_0^2+a_1^2)x_1,
   \\
   f^2(x_1) &= a_0^2+a_1^2+2a_0a_1x_1,
\end{align*}
donde hemos usado, de nuevo, que $x_1^2=1$. Calculando la esperanza de cada una de estas funciones obtenemos
\begin{align*}
   \EE[f^4] &= a_0^4+a_1^4+6a_0^2a_1^2,
   \\
   \EE[f^2] &= a_0^2+a_1^2,
\end{align*}
puesto que $\EE[x_1]=0$ y $\EE[1]=1$. Si $a_1=0$, es decir, $f$ es un polinomio de grado $d=0<1=n$, entonces $\EE[f^4]=\EE[f^2]^2$ y la desigualdad del lema es válida porque $9^0=1$. Finalmente, si $a_1\neq 0$, en el caso cuando $f$ es un polinomio de grado $d=1=n$, entonces
\[
a_0^4+a_1^4+6a_0^2a_1^2 \leq 3\big(a_0^4+a_1^4+2a_0^2a_1^2\big)=3\EE[f^2]^2,
\]
y la desigualdad del lema es válida porque $3<9$.

Supongamos que el lema es cierto en dimensión $n-1$ y queremos probarlo en dimensión $n$.
Primero observamos que podemos reescribir $f$, en términos de su derivada $g=D_nf$ respecto de la $n$-ésima variable, como $f=x_ng + h$ donde $h=f-x_ng$. Notamos que la función $g$ depende de las primeras $n-1$ variables y tiene grado $d-1$, y que la función $h$ depende de las primeras $n-1$ variables y tiene grado menor o igual que $d$.

Basándonos en estas observaciones concluimos que
\begin{align*}
   \|f\|_4^4 &= \EE[f^4] = \EE[(x_ng+h)^4]
   \\
   &= \EE[x_n^4g^4] + 4 \EE[x_n^3g^3h] + 6 \EE[x_n^2g^2h^2]
   + 4\EE[x_ngh^3] + \EE[h^4].
\end{align*}
Recordando que $\EE[x_n^3]=\EE[x_n]=0$, $\EE[x_n^2]=\EE[x_n^4]=1$, y que $x_n$ y $g^kh^j$ son variables independientes, puesto que $g$ y $h$ dependen solo de las primeras $n-1$ variables, obtenemos lo siguiente:
\[
   \|f\|_4^4 = \EE[g^4] + \EE[h^4] + 6\EE[g^2h^2]
   \leq \EE[g^4] + \EE[h^4] + 6 \EE[g^4]^{1/2} \EE[h^4]^{1/2},
\]
donde hemos usado la desigualdad de Cauchy-Schwarz. Como $g$ y $h$ dependen de $n-1$ variables y son respectivamente polinomios de grado $d-1$ y $d$, podemos aplicar la hipótesis inductiva para concluir que
\begin{align*}
   \|f\|_4^4 &\leq 9^{d-1} \EE[g^2]^2 + 9^{d} \EE[h^2]^2 + 6 \cdot 9^{\frac{d-1}{2}}
   9^{\frac{d}{2}} \EE[g^2] \EE[h^2]
   \\
   &\leq 9^{d} \left( \EE[g^2]^2 + \EE[h^2]^2 + 2 \EE[g^2] \EE[h^2] \right)
   \\
   &= 9^d \left( \EE[g^2] + \EE[h^2] \right)^2
   \\
   &= 9^d \EE[f^2]^2.
\end{align*}
En la última igualdad hemos usado una vez más la independencia de las funciones~$g$ y $h$ respecto de la función $x_n$, y que $\EE[x_n]=0$, y por tanto $\EE[g^2] +\EE[h^2]= \EE[x_n^2g^2+h^2]=\EE[(x_ng+h)^2]$.
\end{proof}

Un corolario del lema de Bonami es el llamado \emph{truco de la norma $1$}.

\begin{corollary}[el truco de la norma $1$]\label{cor:1norm-trick}
Sea $f:\{-1,1\}^n\to \R$ un polinomio de grado $d$. Entonces $\|f\|_2\leq 3^{d}\|f\|_1$.
\end{corollary}

\begin{proof} Aplicando la desigualdad de Hölder con los exponentes duales $p=3$ y $p'=3/2$ obtenemos
\[
\|f\|_2^2 = \mathbb{E}\big[ |f(x)|^{4/3}|f(x)|^{2/3}\big] \leq \|f\|_4^{4/3}\|f\|_1^{2/3}.
\]
Por el lema de Bonami (lema~\ref{lem:Bonami's}) tenemos que $\|f\|_4\leq 3^{d/2}\|f\|_2$ y concluimos que
\[
\|f\|_2^2\leq 3^{2d/3}\|f\|_2^{4/3}\|f\|_1^{2/3}.
\]
Reordenando obtenemos el resultado deseado.
\end{proof}

Recientemente los autores de \cite{EKL23} encontraron un análogo al lema de Bonami para funciones \emph{globales} en el espacio de aplicaciones lineales entre dos espacios vectoriales de dimensión finita sobre cuerpos finitos. En este artículo los autores se refieren a la desigualdad hipercontractiva de Bonami-Beckner-Gross, en ese orden, que discutiremos a continuación.

\subsection{La hipercontractividad}\label{sec:hypercontractitvity}

La desigualdad hipercontractiva es un resultado fundamental en análisis, que tiene muchas aplicaciones en matemáticas discretas, teoría de la computación y combinatoria, entre otras disciplinas. Una buena referencia es \cite[capítulo 9 y sección~10.1]{O'Do21} y, para los interesados en el análisis matemático, la tesis doctoral de Bonami \cite{Bon70}, y los artículos de Beckner y Gross \cite{Bec75,Gro75}.

El \emph{operador de ruido} $T_{\rho}$, con $\rho \in \R$, corresponde al siguiente multiplicador de Fourier, definido para todas las funciones $f:\{-1,1\}^n\to\R$ como sigue:
\begin{equation}\label{Trho}
(T_{\rho}f)(x)=\sum_{S\subset [n]} \rho^{|S|} \widehat{f}(S) \chi_S(x).
\end{equation}
Este operador tiene una interesante interpretación probabilística para $|\rho| \leq 1$, que justifica su nombre, y que se puede consultar en \cite[sección 2.4]{O'Do21}. En el contexto de la teoría de la elección social, básicamente nos indica cómo un ruido aleatorio en la recogida de datos de la votación afecta el resultado final.

Observamos que por la unicidad de los coeficientes de Fourier, $\widehat{T_{\rho}f}(S)=\rho^{|S|} \widehat{f}(S)$. En otras palabras, el coeficiente de Fourier de $T_{\rho}f$ correspondiente al subconjunto $S\subset [n]$ es el coeficiente de Fourier de $f$ correspondiente a $S$ multiplicado por un factor que es exponencial en $|S|$, de ahí que de $T_{\rho}$ se diga que
es un \emph{multiplicador de Fourier}. En el contexto del análisis de Fourier en grupos, hay una noción natural de convolución que interactúa con la transformada de Fourier de la manera esperada, es decir, la transformada de Fourier de la convolución de dos funciones es el producto de las transformadas de Fourier de las funciones convolucionadas. Teniendo esto en cuenta, el operador de ruido $T_{\rho}$ es la convolución con una función $R$ cuyos coeficientes de Fourier son $\widehat{R}(S)=\rho^{|S|}$, y esa función es una versión finita del \emph{producto de Riesz} mencionado en el título de la charla que impartió Aline Bonami en Barcelona, definido para cada $x\in\{-1,1\}^n$ como
\[
R(x):=\prod_{i=1}^n (1+\rho x_i)=\sum_{S\subset [n]} \rho^{|S|} \chi_S.
\]

La primera observación es que, gracias a la identidad de Parseval~\eqref{Parseval}, el operador~$T_{\rho}$ es una contracción en $L^2(\{-1,1\}^n)$ cuando $|\rho|\leq 1$:
\[
\|T_{\rho}f \|_2 \leq \|f\|_2.
\]
El teorema de hipercontractividad de Bonami trata sobre la acotación del operador de ruido~\eqref{Trho} de $L^p(\{-1,1\}^n)$ en $L^q(\{-1,1\}^n)$ bajo algunas restricciones en $\rho$, dados $q\geq p\geq 1$.

\begin{theorem}[teorema de hipercontractividad~\cite{Bon70}]
\label{thm:hypercontractivity}
Sean $f:\{-1,1\}^n\to \R$ y $1\leq p \leq q \leq \infty$. Si $\rho^2 \leq \frac{p-1}{q-1}$, entonces
\begin{equation}\label{desigualdad-hiper}
\|T_{\rho} f\|_q\leq \|f\|_p.
\end{equation}
\end{theorem}

Como consecuencia inmediata del teorema tenemos el siguiente resultado en el caso $q=p$:

\begin{corollary}
Si $|\rho|\leq 1$ y $p\geq 1$, entonces $T_{\rho}$ es una contracción en $L^p(\{-1,1\}^n)$, es decir
\[
\|T_{\rho}f \|_p \leq \|f\|_p.
\]
\end{corollary}

El teorema de hipercontractividad nos dice que el operador $T_\rho$, además de ser una contracción de $L^p(\{-1,1\}^n)$ en $L^p(\{-1,1\}^n)$, es una contracción de $L^p(\{-1,1\}^n)$ en $L^q(\{-1,1\}^n)$ para un rango de $q\geq p$. Para ser más precisos, dados $p\geq 1$ y $|\rho |\leq 1$, se cumple la desigualdad $\|T_{\rho}f\|_q\leq \|f\|_p$ para $p\leq q \leq p + \left( \frac{1}{\rho^2} - 1\right) (p-1)$. Es por esto que se habla de hipercontractividad en este contexto.
Observemos que, en el caso $q\leq p$, la desigualdad~\eqref{desigualdad-hiper} vale para todo $|\rho|\leq 1$, y se deduce inmediatamente de la desigualdad de Hölder y del caso contractivo, puesto que
\[
\|T_{\rho}f\|_q\leq \|T_{\rho}f\|_p\leq \|f\|_p.
\]

Otro caso particular de este teorema es el teorema de hipercontractividad-$(4,2)$, que no es más que el teorema de hipercontractividad aplicado a $q=4$, $p=2$ y $\rho=1/\sqrt{3}$. Observamos que $\left( 1/\sqrt{3}\right)^2 = \frac{2-1}{4-1}$, por lo que los parámetros caen en el rango permitido por el teorema~\ref{thm:hypercontractivity}.

\begin{corollary}[teorema de hipercontractividad-$(4,2)$]
\label{cor:(4,2)-hypercontractivity}
Para toda $f:\{-1,1\}^n\to \R$
\[
\| T_{1/\sqrt{3}}f\|_4\leq \|f\|_2.
\]
\end{corollary}

Como veremos a continuación, el teorema de hipercontractividad-$(4,2)$ es equivalente al lema de Bonami. Se deduce inmediatamente de la definición del operador de ruido~\eqref{Trho} que la inversa de $T_{\rho}$ es $T_{1/\rho}$ cuando $\rho\neq 0$ y que si $f$ tiene grado $d$ también lo tendrá $T_{\rho}f$. Ahora bien, asumiendo el teorema de hipercontractividad-$(4,2)$, tendremos que
\[
\|f\|_4=\|T_{1/\sqrt{3}}(T_{\sqrt{3}}f)\|_4\leq \|T_{\sqrt{3}}f\|_2.
\]
Pero, por la identidad de Parseval~\eqref{Parseval}, y como $f$ tiene grado $d$, se obtiene
\[
   \|T_{\sqrt{3}}f\|_2^2
   = \sum_{S\subset [n] \colon |S|\leq d} (\sqrt{3})^{2|S|}|\widehat{f}(S)|^2
   \leq (\sqrt{3})^{2d}\|f\|_2^2,
\]
y tomando la raíz cuadrada en ambos lados de la desigualdad obtenemos la conclusión del lema de Bonami. La implicación contraria, es decir, que el lema de Bonami implica la hipercontractividad-$(4,2)$, es más técnica y puede consultarse en \cite[teorema 9.21 y ejercicio 9.6]{O'Do21}.

El operador de ruido $T_{\rho}$ es autoadjunto, con ello queremos decir que $\langle T_{\rho}f,g\rangle=\langle f,T_{\rho}g\rangle$. Podemos usar esta propiedad para deducir del teorema de hi\-per\-con\-trac\-ti\-vi\-dad-$(4, 2)$ (corolario~\ref{cor:(4,2)-hypercontractivity}) otro caso particular del teorema de hipercontractividad (teorema~\ref{thm:hypercontractivity}), concretamente el correspondiente a los parámetros $p=4/3$, $q=2$ y $\rho=1/\sqrt{3}$, también en el rango permitido por el teorema~\ref{thm:hypercontractivity}. Presentamos este corolario y su demostración a continuación.

\begin{corollary}[teorema de hipercontractividad-$(2,4/3)$]
Para toda función $f:\{-1,1\}^n\to \R$ tenemos que
\[
\| T_{1/\sqrt{3}}f\|_2\leq \|f\|_{4/3}.
\]
\end{corollary}

\begin{proof}
Sea $T=T_{1/\sqrt{3}}$. Si $\|Tf\|_2=0$ no hay nada que probar, la desigualdad es automáticamente cierta. Supongamos que $\|Tf\|_2>0$;
usando que $T$ es autoadjunto, la desigualdad de Hölder para los exponentes duales $p=4/3$ y $p'=4$, y el teorema de hipercontractividad-$(4,2)$, obtenemos que
\[
   \|Tf\|_2^2 = \langle Tf,Tf\rangle = \langle f,T(Tf)\rangle
     \leq \|f\|_{4/3}\|T(Tf)\|_4 \leq \|f\|_{4/3}\|Tf\|_2.
\]
Finalmente, dividiendo por $\|Tf\|_2$, obtenemos la desigualdad deseada.
\end{proof}

El estudio de las desigualdades hipercontractivas en el mundo de las funciones booleanas sigue despertando gran interés. En artículos muy recientes como \cite{KLLM24}, los autores discuten el caso de medidas \emph{$p$-sesgadas}, en contraposición con la medida uniforme (\emph{$1/2$-sesgada}) en el cubo de Hamming en el caso clásico.

\section{La teoría de la elección social y Bonami se encuentran}
\label{sec:social-choice-meets-Bonami}

En esta sección, primero mostraremos que las funciones booleanas no pueden tener todos sus coeficientes de Fourier concentrados en todos los monomios lineales de su desarrollo de Fourier~\eqref{desarrollo_fourier}, de hecho una función afín es o bien una constante, o bien una dictadura o su negativo. A continuación probaremos los teoremas de Friedgut-Kalai-Naor (FKN) y de Kahn-Kalai-Linial (KKL) como consecuencia del lema de Bonami. El teorema de FKN muestra que una función aproximadamente lineal debe estar cerca de una dictadura o su negativo. El teorema de KKL nos permitirá demostrar que, en un sistema de votación monótono con $n$ votantes, el voto a favor de un candidato por parte de una fracción de los votantes (del orden de $1/\log n$) determina, con un $99\%$ de probabilidad, que su candidato sea el elegido.

\subsection{Funciones afines}\label{sec:singletons}

En esta breve sección mostraremos que las funciones booleanas afines son una constante, una dictadura, o el negativo de una dictadura. Una función booleana $f:\{-1,1\}^n\to\{-1,1\}$ es afín si es de la forma $f(x)=a_0+a_1 x_1+\dots+a_n x_n$.

\begin{proposition}\label{prop:boolean-functions-cant-concentrate-singletons}
Si la función booleana $f:\{-1,1\}^n\to \{-1,1\}$ es tal que $\widehat{f}(S)=0$ para todo subconjunto $S\subset [n]$ con $|S|\geq 2$, entonces $f$ es una constante o $f(x)=\pm x_i$ para algún $i\in [n]$.
\end{proposition}

\begin{proof}
Por hipótesis, $f(x)= a_0+a_1x_1+a_2x_2+\dots a_nx_n$. Si la función~$f$ no es constante, entonces existe un $i\in [n]$ tal que $a_i\neq 0$. Pero $D_if(x) = a_i\neq 0$; por tanto $a_i=\pm 1$, puesto que, por ser $f$ una función booleana, $D_if(x)\in \{-1,0,1\}$.
Por la identidad de Parseval~\eqref{Parseval}, $1=\sum_{j=0}^na_j^2=1 +\sum_{j\neq i} a_j^2$; por tanto $a_j=0$ para todo $j\neq i$. De manera que la función $f$ es
$f(x)=a_ix_i$ con $a_i=\pm 1$ y, por tanto, $f$ es la $i$-ésima dictadura $d_i$, o su negativo $-d_i$.
\end{proof}

\subsection{El teorema de Friedgut-Kalai-Naor (FKN)}

Podemos usar el lema de Bonami para probar el teorema de Friedgut, Kalai y Naor~\cite{FKN02}. Este teorema dice que una función booleana que es casi lineal debe estar muy cerca de ser una dictadura o su negativo, es decir, debe estar muy cerca de un carácter lineal o su negativo.

Ser casi lineal significa que la mayoría de la masa espectral reside en los monomios lineales del desarrollo de Fourier de la función $f$. Sea
\[
   \WW^1(f) = \sum_{S\subset [n] \colon |S|=1} \widehat{f}(S)^2
\]
la masa espectral de la parte lineal de $f$. Nótese que, por la identidad de Parseval~\eqref{Parseval}, se tiene que $\WW^1(f)\leq 1$ para toda función booleana. Así que ser casi lineal quiere decir que $\WW^1(f)$ es casi $1$. Cuando $\WW^1(f)=1$, la función $f$ es lineal, es decir, es de la forma $a_1x_1+\dots +a_nx_n$, y en ese caso, como hemos probado en la proposición~\ref{prop:boolean-functions-cant-concentrate-singletons}, $f(x)=\pm x_i$ para algún $i\in [n]$.

\begin{theorem}[FKN]\label{thm:FKN}
Existe una constante absoluta $C>0$ tal que para toda función booleana definida en $\{-1,1\}^n$ se tiene que, para algún $i\in[n]$,
\[
\|f-a_i x_i\|_2^2\leq C(1-\WW^1(f)),
\]
donde $a_i = \widehat{f}(\{i\})$.
\end{theorem}

\begin{proof} Sea $\delta = 1 - \WW^1(f)$. Por la identidad de Parseval \eqref{Parseval}, se obtiene
\[
\|f-a_i x_i\|_2^2 = \sum_{S\subset [n] \colon S \neq \{i\} } \widehat{f}(S)^2 \leq 1,
\]
para cualquier $i\in [n]$. De modo que solo nos interesa el caso de $\delta$ pequeña. Seguiremos el argumento en~\cite[lecture 6]{Min24}.

Para la función booleana
\[
   f(x)=a_0+a_1x_1+\dots + a_nx_n + a_{12}x_1x_2+\dots + a_{1\dots n}x_1\dots x_n,
\]
por hipótesis, tenemos que $\WW^1(f)= a_1^2 +\dots + a_n^2=1-\delta \leq 1$.
Primero afirmamos~que
\begin{equation}\label{claim1}
a_1^4+a_2^4+\dots + a_n^4 \geq 1-C\delta.
\end{equation}
Si la desigualdad \eqref{claim1} es cierta, entonces
\[
   1-C\delta \leq \sum_{i=1}^n a_i^4
   \leq \big(\max_{i\in [n]} a_i^2\big) \sum_{i=1}^n a_i^2
   \leq \max_{i\in [n]} a_i^2,
\]
y por tanto $\max\limits_{i\in [n]} |a_i|^2=|a_{i_0}|^2 \geq 1-C\delta$ para algún $i_0\in [n]$. De aquí se deduce el teorema de FKN, gracias a la identidad de Parseval~\eqref{Parseval}:
\[
   \|f-a_{i_{0}}\chi_{\{i_0\}}\|^2_2=1-|a_{i_0}|^2\leq C\delta.
\]
Para probar la desigualdad~\eqref{claim1}, observemos que
\[
   a_1^4 + \dots + a_n^4 
   = \big(a_1^2+\dots +a_n^2\big)^2 - \sum_{i\neq j} a_i^2a_j^2
   = (1-\delta)^2-\sum_{i\neq j} a_i^2a_j^2.
\]
Como $(1-\delta)^2\geq 1-2\delta$, la desigualdad~\eqref{claim1} estaría probada si pudiéramos demostrar~que
\[
   \sum_{i\neq j} a_i^2a_j^2\leq C\delta.
\]
Consideremos el polinomio de segundo grado
\[
   h(x) = \sum_{i\neq j} a_i a_j x_i x_j
   = \big(a_1x_1+\dots +a_nx_n\big)^2-\big(a_1^2+\dots a_n^2\big).
\]

Primero, observemos que, gracias a la identidad de Parseval~\eqref{Parseval}, se tiene $\|h\|_2^2 = \sum_{i\neq j} a_i^2 a_j^2$. Segundo, notemos que gracias al truco de la norma $1$ (corolario~\ref{cor:1norm-trick}), $\|h\|_2\leq 9 \|h\|_1$. Por lo tanto, bastará demostrar que $\|h\|_1\leq C\sqrt{\delta}$. Por la desigualdad triangular tenemos que
\begin{align*}
   \|h\|_1 &= \big\|
   \big( a_1 x_1 + \dots + a_n x_n \big)^2
   - \big(a_1^2 + \dots + a_n^2 \big)
   \big\|_1
   \\
   &\leq \big\| \big( a_1 x_1 + \dots + a_n x_n \big)^2
   - f(x)^2 \big\|_1
   + \big\| 1 - \big(a_1^2 + \dots + a_n^2 \big) \big\|_1
   \\
   &= \big\| \big(a_1 x_1 + \dots + a_n x_n \big)^2
   - f(x)^2 \big\|_1 + \delta.
\end{align*}
Nótese que en la desigualdad hemos usado por primera vez que $f$ es una función booleana ($f^2(x)=1$). Para estimar el primer término, denotemos la parte lineal $a_1x_1+\dots +a_n x_n$ por $\ell(x)$. Entonces por la desigualdad de Cauchy-Schwarz tenemos~que
\begin{align*}
   \big\| \big(a_1 x_1 + \dots + a_n x_n \big)^2
   - f(x)^2 \big\|_1
   &= \|\ell^2-f^2\|_1
   = \|(\ell-f)(\ell + f)\|_1
   \\
   &\leq \| \ell-f \|_2 \, \| \ell+f \|_2.
\end{align*}
Ahora, $f - \ell = \sum_{S\subset [n] \colon |S|\neq 1} a_S \chi_S$; por lo tanto
\[
   \|\ell-f\|_2^2 = \sum_{S\subset [n] \colon |S|\neq 1} a_S^2 = \delta,
\]
y por la desigualdad triangular,
\[
   \|\ell+f\|_2
   \leq \|\ell\|_2 + \|f\|_2
   = \Big(\sum_{i=1}^n|a_i|^2\Big)^{1/2} +1
   \leq 2.
\]
De modo que
\[
   \|h\|_1 \leq 2 \sqrt{\delta} + \delta \leq 3\sqrt{\delta}.
   \qedhere
\]
\end{proof}

Nuestras estimaciones muestran que en el teorema~\ref{thm:FKN} se puede tomar $C=2+3^6=731$.

\subsection{El teorema de Kahn-Kalai-Linial (KKL)}

En esta sección, $\log$ se refiere al logaritmo neperiano. El teorema de Kahn-Kalai-Linial está enunciado y probado en \cite{KKL88}, y dice que para funciones booleanas que representen el resultado de una elección con $n$ votantes, existe un votante cuya influencia es de orden mayor que $\log n/n$ veces la varianza. Fue en este influyente artículo donde se estableció la conexión entre las funciones booleanas y las desigualdades de Beckner \cite{Bec75} en análisis de Fourier. Los autores, en aquel momento, ignoraban que Aline Bonami había publicado en francés \cite{Bon68,Bon70}, varios años antes, las desigualdades hipercontractivas que ellos necesitaban.

Observemos que gracias a la desigualdad de Poincaré (proposición~\ref{thm:Poincare's Inequality}) y a la definición de influencia tenemos que
\[
   \Var(f)\leq \II(f)\leq n\max\{\II_i(f): i\in [n]\}.
\]
Por tanto, existe un votante $i\in [n]$ cuya influencia satisface:
\[
   \II_i(f)\geq \frac{\Var(f)}{n}.
\]
El teorema de KKL es una mejora espectacular sobre esta estimación.

\begin{theorem}[KKL]\label{thm:KKL}
Existe una constante absoluta $c >0$ tal que para todas las funciones booleanas $f:\{-1,1\}^n\to \{-1,1\}$, existe algún $i\in [n]$ tal que
\[
\II_i(f)\geq c \frac{\log n}{n} \Var(f).
\]
\end{theorem}

Después de probar el teorema de KKL, presentaremos en la sección~\ref{sec:problemaBen-Or_Linial} una de sus consecuencias, la resolución del problema de Ben-Or y Linial.

Para probar el teorema de KKL (teorema~\ref{thm:KKL}) necesitaremos el siguiente lema que es consecuencia del lema de Bonami. Para $f: \{-1,1\}^n\to\R$, definimos $f^{\leq d}$ como el polinomio resultante de truncar el desarrollo de Fourier de $f$ a nivel $d$, es decir
\[
f^{\leq d}(x) = \sum_{S\subset [n], |S|\leq d} \widehat{f}(S) \chi_S(x).
\]

\begin{lemma} \label{lem:auxiliary}
Sea $f: \{-1,1\}^n\to\R$. Entonces $\|f^{\leq d}\|_2^2\leq \sqrt{3}^d \|f\|_2\|f\|_{4/3}$.
\end{lemma}

\begin{proof} Por la ortogonalidad de los caracteres, la desigualdad de Hölder con exponentes duales $p=4$ y $q=4/3$, y el lema de Bonami (lema~\ref{lem:Bonami's}) tenemos que
\[
   \|f^{\leq d}\|_2^2
   = \langle f^{\leq d},f^{\leq d}\rangle
   = \langle f^{\leq d},f\rangle
   \leq \|f^{\leq d}\|_4\|f\|_{4/3}
   \leq \sqrt{3}^d \|f\|_2\|f\|_{4/3}.
   \qedhere
\]
\end{proof}

\begin{proof}[Demostración del Teorema de KKL]
Sea
\[
\theta = \frac{n \max\{\II_i(f): i\in[n]\} }{\Var(f)}.
\]
Por la desigualdad de Poincaré (proposición~\ref{thm:Poincare's Inequality}) se cumple
\[
   \theta \geq \frac{\II(f)}{\Var(f)}\geq 1.
\]
Demostrar el teorema es equivalente a demostrar que $\theta \geq c \log(n)$.

Aplicaremos el lema~\ref{lem:auxiliary} a las funciones
\[
 g_i:= D_if: \{-1,1\}^{n-1}\to \{-1,0,1\}
\]
y aprovecharemos la relación entre la influencia y las derivadas dada en la fórmula~\eqref{eqn:influencia-derivada}. En concreto, la relación
\[
\II_i(f) = \EE[|D_if|^2]=\Prob\big\{|g_i|= 1 \big\}.
\]

Por el lema~\ref{lem:auxiliary} tenemos que
\[
    \|g_i^{\leq d} \|_2^2
    \leq \sqrt{3}^d \|g_i\|_2\|g_i\|_{4/3}
    \leq e^d\big(\II_i(f)\big)^{1/2}\big(\II_i(f)\big)^{3/4}
    = e^d\big(\II_i(f)\big)^{5/4},
\]
porque
\[
   \|g_i\|_p^p=\EE[|g_i|^p]=\Prob \big\{ |g_i|= 1\big\} = \II_i(f).
\]
Si sumamos sobre $i\in [n]$, obtenemos que
\[
\sum_{i=1}^n\|g_i^{\leq d} \|_2^2 \leq e^d \sum_{i=1}^n \big(\II_i(f)\big)^{5/4}
\leq e^d \big(\max_{i\in [n]}\II_i(f)\big)^{1/4} \II(f).
\]
Pero $g_i=\sum_{S\subset [n] \colon i\in S}\widehat{f}(S)\chi_{S\setminus \{i\}}$, entonces $\widehat{g}(S\setminus \{i\})=\widehat{f}(S)$ cuando $i\in S\subset [n]$ y por tanto
\[
   \sum_{i=1}^n\|g_i^{\leq d} \|_2^2
   = \sum_{i=1}^n \sum_{\substack{S\subset [n] \colon i\in S,\\ |S|\leq d+1}} |\widehat{f}(S)|^2
   = \sum_{S\subset [n] \colon |S|\leq d+1} |S| \, |\widehat{f}(S)|^2.
\]
Así que tenemos
\begin{equation}\label{eq:(*)}
   \sum_{S\subset [n] \colon 0<|S|\leq d} |\widehat{f}(S)|^2
   \leq \sum_{S\subset [n] \colon |S|\leq d+1} |S| \, |\widehat{f}(S)|^2
   \leq e^d \II(f) \big(\max_{i\in [n]} \II_i(f)\big)^{1/4}.
\end{equation}

Ahora escogemos $d=\lceil 2 \II(f)/\Var(f) \rceil$, donde $\lceil x \rceil$ representa el menor entero mayor o igual que $x$. Como ya sabemos (por la desigualdad de Poincaré, proposición~\ref{thm:Poincare's Inequality}) que $\II(f)/\Var(f)\geq 1$, el entero $d$ cumple
\[
1 + 2 \II(f)/\Var(f) \geq d \geq 2 \II(f)/\Var(f) \geq 2.
\]
Esta elección garantiza que
\begin{align*}
   \sum_{S\subset [n] \colon |S|> d} |\widehat{f}(S)|^2
   &\leq \frac{1}{d} \sum_{S\subset [n] \colon |S|>d} |\widehat{f}(S)|^2 |S|
   \\
   &\leq \frac{1}{d} \II(f) = \frac{1}{2 \II(f)/\Var(f)} \II(f)
   = \frac{1}{2} \Var(f).
\end{align*}
Por lo tanto,
\[
   \sum_{S\subset [n] \colon 0<|S|\leq d} |\widehat{f}(S)|^2
   = \Var(f)- \sum_{S\subset [n] \colon |S|> d}|\widehat{f}(S)|^2
   \geq \frac{1}{2} \Var(f)
\]
y por la desigualdad \eqref{eq:(*)} y la cota superior en $d$ tenemos que
\[
   \big(\max_{i\in [n]}\II_i(f)\big)^{1/4}
   \geq \frac{1}{2e} \frac{\Var(f)}{\II(f)} e^{-2\frac{\II(f)}{\Var(f)}}
   \geq \frac{1}{2e} e^{-3\frac{\II(f)}{\Var(f)}}.
\]
La última desigualdad es consecuencia de que $x\leq e^x$ para todo $x\in \R$. Elevando a la cuarta potencia, concluimos que existen constantes $c_1,c_2>0$ tales que
\begin{equation}\label{eq:max_Ii}
\max\{\II_i(f): i\in[n]\} \geq c_1e^{-c_2\frac{\II(f)}{\Var(f)}}.
\end{equation}
Utilizando que $n \max\{\II_i(f): i\in[n]\} \geq \II(f)$ y la definición de $\theta$, obtenemos
\[
\theta \Var(f) \geq c_1 n e^{-c_2 \theta}.
\]
Pero la varianza de $f$ es menor que $1$, así que
\[
\theta \geq c_1 n e^{-c_2 \theta}.
\]
Esta relación la reescribimos en la forma
\begin{equation}\label{eq:theta}
\theta e^{c_2 \theta} \geq c_1 n .
\end{equation}
Ahora, es un breve cálculo mostrar que la desigualdad \eqref{eq:theta} para un valor $\theta \geq 1$ implica que
\[
\theta \geq c \log n,
\]
para una constante $c = c (c_1, c_2)$.
\end{proof}

Si seguimos las constantes, concluimos que podemos tomar
\[
   c \leq (\log(c_1^{-1})+c_2+1)^{-1}
\]
y las constantes $c_1^{-1}=(2e)^{4}$, $c_2=12$, con lo cual podemos estimar $c=1/20= 0.05$. Una estimación más fina de la constante da un valor de $1/2$ para $n$ grande. Los detalles se pueden encontrar en \cite[ejercicio 9.30]{O'Do21}.

Ben-Or y Linial \cite{BL85} demostraron cómo elegir una partición $\mathcal{P}$ de $[n]$ de manera~que
\[
   \max\{\II_i(\Tribus_{\mathcal{P}}):i\in [n]\}
   \simeq \frac{\log n}{n} \Var\big[\Tribus_{\mathcal{P}}\big]
\]
para la función $\Tribus_{\mathcal{P}}$ introducida en el ejemplo~\ref{eg:tribes}. La partición es uniforme y consiste en $w$ bloques de $s$ elementos, donde $s=s(w)$ ha sido definido como el mayor número entero tal que $1-(1-2^{-w})^s \leq 1/2$, y la dimensión es $n=sw$. Esta selección de parámetros hace que la función casi no tenga sesgo, ya que la probabilidad de que la función sea $-1$ está lo más cerca de $1/2$ posible.
En este ejemplo, cada votante tiene poca influencia, de hecho
\[
\II_i(\Tribus_{\mathcal{P}}) = \frac{\log n}{n} \big(1\pm o(1)\big)
\]
para cada $i\in [n]$. Por otro lado, $\EE[\Tribus_{\mathcal{P}}]=o(1)$, de modo que $\Var[\Tribus_{\mathcal{P}}] = 1 - o(1)$, y este ejemplo muestra que la cota inferior en el teorema KKL es óptima. Se puede consultar la referencia~\cite[sección 4.2]{O'Do21} para más detalles.

\subsection{El efecto del voto de unos pocos}\label{sec:problemaBen-Or_Linial}

En esta sección estudiamos cómo, asegurando el voto de una fracción relativamente pequeña de la población a favor de un candidato, en un sistema de votación representado por una funcion booleana monótona, se puede sesgar el resultado de la elección a favor de ese candidato. Para asegurar esos votos se puede intervenir en debates, distribuir publicidad positiva, o alternativamente saturar las redes de bulos o incluso chantajear o sobornar a una fracción de la población.

El problema de Ben-Or y Linial dice que para una función booleana que es monótona y no es constante, que representa el resultado de una elección con $n$ votantes, existe una coalición del orden de $n/ \log n$ votantes cuyos votos pueden controlar el resultado de la elección con probabilidad muy cerca de $1$. Para ser más precisos, podemos determinar cuántos y quiénes necesitan votar a favor de un candidato para asegurar el resultado de una elección con un $99\%$ de certidumbre, sabiendo que inicialmente el candidato al que se quiere favorecer tiene por lo menos un $0.5\%$ de probabilidad de ganar la elección. El lector interesado puede cambiar los porcentajes y repetir los cálculos, obteniendo diferentes constantes de comparación. Formalizamos este resultado en un teorema.

\begin{theorem}[el problema de Ben-Or y Linial]\label{thm:BenOr-Linial}
Sea $f:\{-1,1\}^n\to\{-1,1\}$ una función booleana monótona tal que $\EE[f]\geq -0.99$. Entonces existe un subconjunto $J\subset [n]$ con $|J|\leq O(n/\log n)$ tal que si los \emph{convencemos para que voten} $1$ forzarán el resultado final a ser $1$ casi siempre, es decir,
\[
\EE[f(x_1,x_2,\dots,x_n)= 1: x_i= 1 \ \forall i\in J] \geq 0.99.
\]
De manera similar, si $\EE[f]\leq 0.99$, entonces existe un subconjunto $J\subset [n]$ con $|J|\leq O(n/\log n)$ tal que si los \emph{convencemos para que voten} $-1$ forzarán el resultado final a ser~$-1$ casi siempre, es decir,
\[
\EE[f(x_1,x_2,\dots,x_n)=-1: x_i=-1 \ \forall i\in J] \leq -0.99.
\]
\end{theorem}

Para una versión más precisa de este resultado ver \cite[ejercicio 9.27]{O'Do21}.

\begin{proof}
Sin perder la generalidad, consideremos el primer caso. Por tanto tenemos una función booleana monótona $f$ definida en $\{-1,1\}^n$ y tal que $\EE[f] \geq -0.99$ (esto ciertamente implica que la función no es constantemente igual a $-1$, más aún, es equivalente a decir que el candidato representado por $1$ tiene al menos un $0.5\%$ de probabilidad de ganar. Vale la pena recordar que si $f$ es una función booleana, $\EE[f] = \Prob\{f=1\} - \Prob\{f=-1\} = 2 \Prob\{f=1\} - 1 \geq -0.99$ por lo que $\Prob\{f=1\}\geq (1-0.99)/2=0.005$).

Seguiremos la siguiente estrategia recursiva para identificar los votantes que han de votar por $1$: siendo $f_0:=f$ la función booleana monótona dada, definida en $\{-1,1\}^n$ y tal que $\EE[f_0]\geq -0.99$, identificamos al votante $i_1$ que tenga la mayor influencia en $f_0$ y nos aseguramos de que vota $1$. Actualicemos la función $f_1:=f_0^{(i_1\to 1)}$, una nueva función booleana monótona definida ahora en $\{-1,1\}^{n-1}$ (todos los votantes originales salvo por el votante $i_1$, cuyo voto es ahora $1$); a continuación identificamos al votante $i_2$ con la mayor influencia en $f_1$. Sea $f_2:=f_1^{(i_2\to1)}=f_0^{(i_1,i_2\to 1)}$ una nueva función booleana monótona definida en $\{-1,1\}^{n-2}$. Continuemos de manera recursiva: después de $k<n$ pasos, $k$ votantes ahora apoyan al candidato $1$, los correspondientes a los índices $i_1, i_2,\dots, i_k$, y tenemos $k+1$ funciones $f_0,f_1,f_2,\dots, f_k$, cada una definida en una variable menos que la anterior, siendo $f_j$ una función de $n-j$ variables. Actualizamos la función a $f_{k+1}=f_{k}^{(i_{k+1}\to 1)}$, una nueva función booleana monótona definida en $\{-1,1\}^{n-(k +1)}$ (todos los votantes originales salvo los $k+1$ votantes cuyos votos son ahora $1$) y repetimos.

Al cabo de $n$ pasos, todos los votantes van a votar por $1$: como la función es monótona y no es idénticamente igual a $-1$, la función no tiene más remedio que ser~$1$. Queremos argumentar a continuación que, cuando $n$ es muy grande, mucho antes de que el proceso finalice, la función será $1$ con una probabilidad del $99\%$, y para ello usaremos el teorema de KKL.

Por la manera de construir las funciones $f_k$, para $k\leq n$ tenemos que
\begin{equation}\label{eq:bribe}
\EE[f_{k+1}] = \EE[f_k] + \max\{\II_i [f_k]: i\in [n-k]\}.
\end{equation}
Esta igualdad muestra que la esperanza de la función $f_k$ aumenta en cada paso, es decir,
\[
\EE[f_{k+1}] \geq \EE[f_k] \text{ para cada } k\leq n-1.
\]
En particular, $\EE[f_k]\geq \EE[f_0]\geq -0.99$. Podemos ver que la igualdad \eqref{eq:bribe} se verifica porque la esperanza de una función booleana monótona $g$ definida en $\{-1,1\}^j$ es el coeficiente de Fourier correspondiente al conjunto vacío y la influencia individual es $\II_i(g)=\widehat{g}(\{i\})$, gracias a la proposición~\ref{prop:influenceMonotone}. Para ser más precisos, si
\[
   g(x) = a_0+a_1x_1+\dots + a_jx_j +\sum_{S\subset [j] \colon |S|>1} a_Sx^S,
\]
entonces $\EE[g]=a_0$. Cuando el votante (digamos el votante $i$) cuya influencia es la mayor vota $1$, y declaramos $g_1:=g^{(i\to1)}$, entonces $\EE[g_1]=a_0+a_{i}=\EE[g]+\II_{i}(g)$, pero por construcción $\II_{i}(g)= \max\{\II_{\ell} (g): \ell\in [j]\}$. Aplicamos este resultado a $g=f_{k}$ y $j=n-k$ y obtenemos \eqref{eq:bribe}.

Para obtener el resultado que buscamos, bastará encontrar $j$ tal que $\EE[f_j]\geq 0.99$ (esto implica que $f_j(x)=1$ con probabilidad $0.995$, por un cálculo similar al hecho en el primer párrafo de la prueba).
Supongamos que después de $k$ pasos $\EE[f_k]< 0.99$; puesto que ya hemos visto que $\EE[f_k]\geq -0.99$, podemos concluir que $\Var [f_k]=1-\EE^2[f_k]=c_0>0$.
Por otro lado, por el teorema~\ref{thm:KKL} tenemos que
\[
   \max\{\II_i(f_k): 1\leq i\leq n-k\}\geq c c_0 \frac{\log (n-k)}{n-k};
\]
lo mismo ocurre en los pasos anteriores. Sea $c_1= c c_0 $; aplicando recursivamente la igualdad~\eqref{eq:bribe} obtenemos que
\begin{align*}
   \EE[f_{k+1}] &= \EE[f_0] + \sum_{j=1}^k \max\{\II_i [f_j]: i\in [n-j]\}
   \\
   &\geq c_1\sum_{j=1}^k \frac{\log (n-j)}{n-j} -0.99
   \geq c_1 k\frac{\log n}{n} - 0.99,
\end{align*}
donde hemos usado que la función $(\log t)/t$ es decreciente para $t>2$. Concluimos que $\EE[f_{k+1}]\geq 0.99$ después de que $k = \lceil 1.98 n / (c_1\log n)\rceil$ votantes voten $1$. En otras palabras, si nos aseguramos de que una coalición de orden $n/\log n$ vote $1$, ¡el resultado de la elección es $1$ con un $99.5\%$ de certidumbre, mejor que lo prometido!
\end{proof}

El lector curioso se preguntará a cuántos votantes ha de convencer, en una sociedad de por ejemplo $n=10^6$ votantes, para que gane su candidato con $99\%$ de certeza. El resultado que acabamos de presentar es útil si el orden de magnitud de~$n$ es mucho más grande que $10^6$, esto es debido principalmente a que las constantes se acumulan y además no son óptimas. Invitamos al lector interesado a explorar con más detalle las estimaciones.

Para más información sobre este problema se pueden consultar las referencias~\cite{BL85} y \cite[proposición 9.27, p. 265--266]{O'Do21}. Añadimos por último que una versión del teorema~\ref{thm:BenOr-Linial}, para funciones booleanas que no son necesariamente monótonas, se puede encontrar en \cite[ejercicio 9.28]{O'Do21}.

\section{Entropía vs. influencia}\label{sec:EntropyInfluence}

En esta sección presentamos las nociones de entropía y mínima entropía para funciones booleanas y la famosa conjetura \emph{entropía de Fourier / influencia} (EFI), planteada por Friedgut y Kalai en 1996 \cite{FK96}. Veremos algunas de las técnicas que se han utilizado para estudiar la conjetura en determinadas clases de funciones, donde o bien se ha demostrado, o se han obtenido resultados parciales. También discutimos la conjetura más débil \emph{mínima entropía de Fourier / influencia} (MEFI) y su conexión con el teorema KKL.

A partir de ahora, $\log$ denota siempre el logaritmo en base $2$. Establecemos además que $0 \log 0=0$.

\subsection{Entropía de Fourier, conjeturas EFI y MEFI }

Dada una distribución de probabilidad $p$ sobre un conjunto $\{z_1,z_2,\dots\}$ numerable, con $p_i=p(z_i)$, la \emph{entropía de Shannon} de $p$ se define como
\[
H(p):=\sum_{i\geq 1} p_i \log \frac{1}{p_i}.
\]
La cantidad $H(p)$ representa una medida de la información transportada por $p$, donde una mayor entropía corresponde a una menor información (más incertidumbre o más falta de información). Un ejemplo básico es el del lanzamiento de una moneda: si $p(\text{cara})=1$ y $p(\text{cruz})=0$, la entropía es $0$, y no hay incertidumbre alguna. Por otro lado, la máxima entropía se alcanza si $p(\text{cara})= p(\text{cruz})= 1/2$, es decir cuando la aleatoriedad es máxima. En general, si el conjunto es finito, la entropía es máxima cuando la distribución $p$ es uniforme, y es 0 cuando $p$ es una distribución puntual.

A continuación definimos la entropía y la mínima entropía de Fourier asociadas a funciones en el cubo de Hamming.

\begin{definition} Dada $f:\{-1,1\}^n \to \R$ con $\|f\|_2=1$, definimos su \emph{entropía de Fourier}\footnote{Obsérvese que coincide con la entropía de Shannon de la distribución de probabilidad $|\widehat{f}(S)|^2$ definida en $S\subset [n]$.} (también denominada \emph{entropía espectral}) $\HH(f)$ como
\[
\HH(f) = \sum_{S\subset [n]} |\widehat{f}(S)|^2\log \frac{1}{|\widehat{f}(S)|^2}.
\]
\end{definition}

Recordemos que si $f$ es una función booleana, $\|f\|_2=1$.

\begin{definition} Dada $f:\{-1,1\}^n \to \R$ con $\|f\|_2=1$, la \emph{entropía mínima de Fourier} $\HH_{\infty}(f)$ se define como
\[
\HH_{\infty}(f) = \min_{S\subset [n]} \log \frac{1}{|\widehat{f}(S)|^2}.
\]
\end{definition}

Intuitivamente,
la entropía de Fourier mide la dispersión de la distribución de Fourier sobre los $2^n$ subconjuntos de $[n]$, y la influencia total mide la concentración de la distribución de Fourier en los coeficientes de nivel alto.
La conjetura EFI nos dice que para funciones booleanas, si la distribución de Fourier es muy dispersa (es decir, funciones con una gran entropía de Fourier) entonces los monomios de alto grado deben tener un peso de Fourier significativo (es decir, su influencia total es grande). Recordemos que la influencia total de una función booleana $f:\{-1,1\}^n\to \{-1,1\}$ viene dada en términos de los coeficientes de Fourier de $f$ por $\II(f)=\sum_{S\subset [n]} |\widehat{f}(S)|^2 |S|$ (ver la proposición~\ref{prop:influence-Fourier}).

\begin{conjecture}[conjetura entropía de Fourier / influencia (EFI) \cite{FK96}]
Existe una constante $c>0$ tal que para toda $f:\{-1,1\}^n\to\{-1,1\}$ se verifica
\[
\HH(f)\leq c \II(f).
\]
\end{conjecture}

Una primera observación es que la entropía y la influencia de $f$ están ambas acotadas por $n$. En efecto, como $\widehat{f}(S)=\frac{1}{2^n}\sum f \chi_S$ y $\sum f \chi_S$ es un número entero, si $\widehat{f}(S)\neq 0$ entonces $|\widehat{f}(S)|\geq 1/2^n$, y por tanto $\HH(f) \leq n$. Por otro lado, como $|S|\leq n$, tenemos que $\II(f)\leq n$.

La conjetura EFI es falsa para funciones con valores reales en el cubo de Hamming. Por ejemplo, concentrando todos los coeficientes de Fourier en el primer nivel, digamos $\widehat{f}(\{x_i\})=1/\sqrt{n}$ para todo $i\in [n]$ y $0$ para el resto de los coeficientes, obtenemos una función $f$ con $\HH(f)=\log(n)$ e $\II(f)=1$. Un contraejemplo de este tipo no es posible si la función es booleana, ya que las únicas funciones booleanas cuyos coeficientes de Fourier se concentran en el primer nivel son de la forma $f(x)=x_i$ para algún $i\in [n]$ (ver la sección~\ref{sec:singletons}).

Se puede comprobar que la conjetura es cierta para las funciones booleanas canónicas analizadas en las secciones anteriores. Es trivial para las funciones paridad. Se puede verificar también, mediante cálculos más o menos simples, para las funciones $\AND$ y $\OR$, la función mayoría (consultar la sección~\ref{sec:canonical-examples}) y la función $\Tribus$ \cite{BL85}. En todos estos casos tenemos la ventaja de que existe una fórmula explícita para los coeficientes de Fourier.

Entre las clases de funciones más generales para las que se ha demostrado la conjetura están las funciones simétricas
\cite{OWZ11}. En el próximo apartado comentaremos algunas de las técnicas que se utilizan en estas demostraciones, así como otros resultados parciales.

\begin{conjecture}[conjetura mínima entropía de Fourier / influencia (MEFI)]
Existe una constante $c>0$ tal que para toda $f:\{-1,1\}^n\to\{-1,1\}$ se verifica
\[
\HH_{\infty}(f)\leq c \II(f).
\]
\end{conjecture}

Es inmediato comprobar que la conjetura MEFI es más débil que la conjetura~EFI.

La condición $\HH_{\infty}(f)\leq c \II(f)$ es equivalente a la existencia de un subconjunto $S\subset [n]$ tal que $|\widehat{f}(S)|^2\geq 2^{-c \II(f)}$. De nuevo, como muestra el ejemplo anterior, la conjetura MEFI no es cierta para funciones con valores reales en el cubo de Hamming.

El teorema de KKL para funciones sin sesgo, es decir funciones con $\EE[f]=0$, se puede deducir de la conjetura MEFI \cite{ACK+20, OWZ11}. Para ver esto, observemos que dado $S\neq \emptyset$, para todo $j \in S$ se tiene $\II_j(f)\geq |\widehat{f}(S)|^2$, por la fórmula de Fourier para la influencia individual~\eqref{eq:influencia_i}, con lo cual
\[
\max\{ \II_i (f): i\in [n]\} \geq \II_j(f)\geq |\widehat{f}(S)|^2.
\]
Como $\widehat{f}(\emptyset)=\EE[f]=0$, y si la conjetura MEFI es válida, esta implica
\[
\max \{ \II_i (f): i\in [n]\} \geq | \widehat{f}(S)|^2\geq e^{-c \II(f)},
\]
que es exactamente la expresión~\eqref{eq:max_Ii} en el teorema de KKL para funciones sin sesgo, puesto que en este caso $\Var[f]=1$.

Podemos deducir el teorema de KKL general (para funciones booleanas con sesgo) de la conjetura EFI. A continuación explicamos el argumento.
Tenemos que
\begin{align*}
   \min_{S\subset [n] \colon S\neq \emptyset}
   \Big(\log\frac{1}{|\widehat{f}(S)|^2}\Big) \Var[f]
   &\leq \sum_{S\subset [n] \colon S\neq \emptyset} |\widehat{f}(S)|^2
   \log\frac{1}{|\widehat{f}(S)|^2}
   \\
   &\leq \sum_{S\subset [n]} |\widehat{f}(S)|^2 \log\frac{1}{|\widehat{f}(S)|^2}
   \\
   &\leq c \II(f),
\end{align*}
donde la última desigualdad vale porque estamos asumiendo que la conjetura EFI es cierta y hemos usado la identidad $\Var[f]=\sum_{S\subset [n], S\neq\emptyset} |\widehat{f}(S)|^2$.

Por tanto,
\[
   \min_{S\subset [n] \colon S\neq \emptyset}
   \log\frac{1}{|\widehat{f}(S)|^2}
   \leq c \frac{\II(f)}{\Var[f]}.
\]
Equivalentemente, para algún $S\neq \emptyset$ tenemos que
\[
|\widehat{f}(S)|^2 \geq 2^{-c\frac{\II(f)}{\Var[f]}}.
\]
El resto del argumento es el mismo que hicimos en el caso $\EE[f]=0$.

La conjetura MEFI se cumple para funciones booleanas monótonas, si no son demasiado sesgadas. Curiosamente, esto se deduce trivialmente del teorema de KKL \cite{OWZ11} como veremos a continuación. Sin embargo, la conjetura EFI sigue abierta para esta clase de funciones.

Supongamos que $f$ es una función monótona booleana con $\Var (f)\geq c_0$ para algún $c_0>0$. Recordemos que por ser $f$ monótona, $\widehat{f}(i)=\II_i(f)$ (ver la proposición~\ref{prop:influenceMonotone}). Entonces, por el teorema de KKL~\eqref{eq:max_Ii}, podemos deducir que
\[
   \max_{S\subset [n]} |\widehat{f}(S)|
   \geq \max_{i\in [n]} \widehat{f}(i)
   = \max \{ \II_i (f): i\in [n]\}
   \geq e^{-c'\II(f)/\Var(f)} \geq e^{-c \II(f)},
\]
donde $c$ es una constante que depende de $c_0$.

Recientemente, en \cite{KKK+20}, los autores han presentado una prueba del teorema de KKL sin usar desigualdades hipercontractivas o las conjeturas EFI/MEFI. Usan \emph{solamente} la desigualdad Log-Sobolev establecida por
Gross \cite{Gro75}, que es equivalente a la desigualdad de hipercontractividad. Gross \cite{Gro75} mostró que la desigualdad Log-Sobolev implica el teorema de hipercontractividad y la implicación inversa puede encontrarse en \cite[ejercicio 10.23]{O'Do21}.

\subsection{La conjetura EFI}

En esta sección analizaremos resultados relacionados con la conjetura EFI. El principal objetivo es mostrar una pequeña parte del amplio abanico de técnicas y puntos de vista diferentes utilizados en los numerosos intentos de demostrar la conjetura.

Dada una $f:\{-1,1\}^n \to \R$, definimos la masa espectral de $f$ en el nivel $k$ como:
\[
\WW^k(f)= \sum_{S\subset [n] \colon |S|=k} |\widehat{f}(S)|^2.
\]
La influencia de $f$ se puede expresar en función de la masa espectral de $f$ en los niveles $k$ como
\[
\II(f) = \sum_{k=0}^n k \WW^k(f).
\]

Si la función $f$ es booleana, por la identidad de Parseval~\eqref{Parseval}, $\sum_{k=1}^n\WW^k(f)=1$. El siguiente resultado, que puede ser interpretado por un analista como una versión tipo \emph{Littlewood-Paley en el contexto de la conjetura EFI}, establece la relación entre la entropía de Shannon de la distribución de probabilidad
$\{ \WW^k (f)\}_k$ y la influencia de $f$, $\II(f)$.

\begin{theorem}[teorema 5 en \cite{OWZ11}]\label{A}\label{thm:OWZ11-Theorem 5}
Sea $f: \{-1,1\}^n\to \{-1,1\}$. Entonces
\[
\sum_{k=0}^n \WW^k(f) \log \frac{1}{\WW^k(f)}\leq 3 \II(f).
\]
\end{theorem}

La pieza clave de la demostración es la denominada \emph{desigualdad isoperimétrica de aristas}\footnote{\emph{Edge-isoperimetric inequality} en la literatura matemática en inglés.} en el cubo de Hamming. Esta desigualdad, muy relevante en el estudio de funciones booleanas, nos dice que si denotamos por $p$ la probabilidad de $f=1$, entonces $\II(f)\geq 2p \log(1/p)$. Observemos que $\II(f)$ representa la proporción de aristas del cubo que unen vértices contiguos en los que $f$ cambia de valor (ver \cite[proposición~2]{OWZ11} y también \cite{FS07}).

El teorema~\ref{A} nos permite abordar una primera aproximación a la conjetura EFI donde, en lugar de obtener una constante independiente de la dimensión $n$, obtenemos $O(\log n)$.

\begin{proposition}\label{logn}
Existe una constante $c>0$ tal que si $f:\{-1,1\}^n\to\{-1,1\}$, entonces $\HH(f)\leq c (\log n) \II(f)$.
\end{proposition}

\begin{proof}
Utilizaremos la notación $W_k$ en lugar de $\WW^k(f)$. La entropía de~$f$ se puede reescribir como
\begin{align*}
   \HH(f) &= \sum_{S\subset [n]} |\widehat{f}(S)|^2 \log \frac{1}{|\widehat{f}(S)|^2}
   \\
   &= \sum_{k=0}^n \bigg(
   \sum_{|S|=k} |\widehat{f}(S)|^2 \log\frac{1}{|\widehat{f}(S)|^2}
   \bigg)
   \\
   &= \sum_{k=0}^n W_k\bigg(
   \sum_{|S|=k} \frac{|\widehat{f}(S)|^2}{W_k} \log\frac{1}{|\widehat{f}(S)|^2}
   \bigg).
\end{align*}

Observemos que para $k=0,1,\dots, n$, $\{ |\widehat{f}(S)|^2/W_k: |S|=k\}$ define una medida de probabilidad cuando $W_k>0$ (si $W_k=0$ lo descartamos y el término no aparece en la suma). Puesto que $\log x$ es una función cóncava, aplicando el resultado análogo a la desigualdad de Jensen, obtenemos
\begin{align*}
   \HH(f) &\leq \sum_{k=0}^n W_k \log\bigg(\frac{1}{W_k}\sum_{|S|=k} 1\bigg)
   \\
   &= \sum_{k=0}^n W_k \log \left( \frac{1}{W_k} \binom{n}{k} \right)
   \\
   &= \sum_{k=0}^n W_k \log \left(\frac{1}{W_k}\right)
   + \sum_{k=0}^n W_k \log \binom{n}{k}.
\end{align*}
Por el teorema~\ref{thm:OWZ11-Theorem 5},
\[
\sum_{k=0}^n W_k \log \left(\frac{1}{W_k}\right)\leq 3 \II(f).
\]
Por otro lado, la conocida desigualdad sobre números combinatorios $\binom{n}{k} \leq \left(\frac{en}{k}\right)^k$ nos permite estimar, para $n\geq 3$,
\[
   \log\binom{n}{k}
   \leq k \log \left(\frac{en}{k}\right)
   \leq k\log e + k\log n
   \leq 2k \log n.
\]
Por tanto,
\[
   \sum_{k=0}^n W_k \log \binom{n}{k}
   \leq 2(\log n) \sum_{k=1}^n k W_k = 2(\log n) \II(f).
   \qedhere
\]
\end{proof}

Si la función booleana $f$ es \emph{simétrica}, es decir, invariante por permutaciones en~$[n]$, se puede obtener una estimación más precisa del término $\sum_{k=0}^n W_k \log \binom{n}{k}$. En \cite[teorema 4]{OWZ11} demuestran que
\[
\sum_{k=0}^n W_k \log \binom{n}{k} < 10 \II(f),
\]
probando así la conjetura EFI para funciones simétricas.

La noción de simetría se puede generalizar a las llamadas \emph{funciones $d$-simétricas}. En este caso, existe una partición $[n] = V_1\cup V_2\cup\dots\cup V_d$ de forma que $f$
es invariante bajo permutaciones de las coordenadas en cualquiera de los $V_i$. Utilizando razonamientos similares al caso simétrico, en \cite[teorema 2]{OWZ11} obtienen para funciones booleanas $d$-simétricas el siguiente resultado parcial:
$\HH(f) \leq C(d) \II(f),$
donde $C(d)=12.04+\log d$.

Un resultado similar al mostrado en la proposición \ref{logn} se verifica también para funciones definidas en el cubo de Hamming, no necesariamente booleanas, siempre y cuando $\|f\|_2 = 1$.

\begin{proposition}\label{prop:logn}
Existe una constante $c>0$ tal que si $f:\{-1,1\}^n\to\mathbb{R}$ y $\|f\|_2 = 1$, entonces $\HH(f)\leq c \big( \II(f)\log n+1\big)$.
\end{proposition}

Las funciones $f:\{-1,1\}^n\to \mathbb{R}$ con $\widehat{f}(\emptyset)^2=1-\epsilon$, $\widehat{f}(\{x_1\})^2=\epsilon$, cero el resto de coeficientes, y $\epsilon \to 0$ muestran que la desigualdad es óptima, en el sentido de que no se puede quitar el $+1$ del lado derecho. La demostración que presentamos a continuación está basada en la teoría matemática de la comunicación de Shannon. Su conexión con la relación entropía-influencia ha sido estudiada en \cite{WWW14}. El concepto clave es el de \emph{código fuente de una variable aleatoria}.

Consideremos una variable aleatoria $X$ que toma valores en un conjunto $\Delta$ con distribución de probabilidad $p$. La entropía (de Shannon) de $X$ es
\[
H(X) = \sum_{x\in\Delta} p(x)\log\frac{1}{p(x)}.
\]

Para un alfabeto finito (conjunto de símbolos) $A$, el conjunto de todas las cadenas finitas formadas por los elementos de $A$ se denota $A^*$. La longitud de una cadena $a\in A^*$ se escribe $|a|$. Un código fuente $C$ para la variable aleatoria $X$ es simplemente una función de $\Delta$ a $A^*$. Un código sin prefijo es aquel en el que $C(x)$ no es un prefijo de $C(y)$ siempre que $x\neq y$.

El teorema de codificación de Shannon \cite{Sha48} examina la relación entre $H(X)$ y la longitud esperada de $C(X)$. El resultado que nos interesa es el límite superior de la entropía.

\begin{theorem}[\cite{Sha48}] Para un código sin prefijo $C$,
\[
H(X)\leq \log (|A|) \EE[|C(X)|].
\]
\end{theorem}

\begin{proof}[Demostración de la proposición~\ref{prop:logn}]
Supongamos ahora que $f$ es una función definida en el cubo de Hamming $\{-1, 1\}^n$ que toma valores reales. Si $\|f\|_2 = 1$, entonces por la identidad de Parseval~\eqref{Parseval} los coeficientes de Fourier $\{|\widehat{f}(S)|^2\}$ definen una medida de probabilidad en los subconjuntos $S$ de $[n]$. Sea $X$ la variable aleatoria que asigna a cada $S\subset [n]$ la probabilidad $|\widehat{f}(S)|^2$. Entonces $H(X)=\HH(f)$ y, por el teorema de Shannon,
\[
\HH(f)\leq \log (|A|) \sum_{S\subset [n]} |\widehat{f}(S)|^2 |C(S)|
\]
para todo código sin prefijo $C$, con alfabeto $A$ definido en los subconjuntos de $[n]$.

A continuación, definimos un código $C$ en los subconjuntos de $[n]$. Cada $i\in [n]$ tiene un desarrollo binario único, $b(i)$, de longitud exactamente\footnote{Los ceros a la izquierda forman parte del desarrollo binario.} $\lceil\log n\rceil$, recordando que $\lceil x\rceil$ denota el menor entero mayor que $x$. Definimos $C(S)$ como la cadena formada al concatenar $b(i)$ para todos los $i\in S$, con $0$ asignado a la cadena vacía. Entonces $|C(S)| = \lceil\log n\rceil |S|$ y tenemos
\[
   \sum_{S\subset [n]} |\widehat{f}(S)|^2 |C(S)|
   = \lceil\log n\rceil \sum_{S\subset [n]} |\widehat{f}(S)|^2 |S|
   = \lceil\log n\rceil \II(f).
\]
Ahora nos gustaría aplicar el resultado de Shannon para obtener el límite logarítmico de $\HH(f)$; sin embargo, nuestro código $C$ no está libre de prefijos: $C(\{1\})$ es un prefijo de $C( \{1,2\})$, por ejemplo. Para tener un código sin prefijos, agregamos el símbolo $!$ a nuestro alfabeto y lo colocamos al final de cada cadena producida por el código. Este código modificado utiliza el alfabeto de tres letras $\{0,1,! \}$, no tiene prefijos y asigna a $S$ una cadena de longitud $\lceil\log n\rceil |S| + 1$. Aplicando el teorema de Shannon, obtenemos
\[
\HH(f)\leq \log 3 \left(\lceil\log n\rceil \II(f) + 1 \right).
\qedhere
\]
\end{proof}

La misma técnica, con un poco de trabajo adicional, permite recuperar el resultado establecido anteriormente para funciones booleanas eliminando el $+1$. Los detalles se pueden encontrar en \cite{WWW14}.

El código fuente anterior no depende de ninguna manera de la naturaleza de $f$, se define puramente en términos de los conjuntos $S$. Si para funciones booleanas pudiéramos definir un código fuente cuya longitud fuera, en promedio, proporcional a la media de $|S|$, es decir,
\[
   \sum_{S\subset [n]} |\widehat{f}(S)|^2 |C(S)|
   \sim \sum_{S\subset [n]} |\widehat{f}(S)|^2 |S|,
\]
entonces obtendríamos inmediatamente la conjetura EFI, puesto que
\[
   \HH(f)
   \leq \log (|A|) \sum_{S\subset [n]} |\widehat{f}(S)|^2 |C(S)|
   \leq K \sum_{S\subset [n]} |\widehat{f}(S)|^2 |S|
   = K \II(f).
\]

La conjetura EFI ha sido verificada para varias familias de funciones booleanas, entre ellas las funciones simétricas o $d$-simétricas con $d$ constante, mencionadas anteriormente, los llamados \emph{árboles de decisión de una sola lectura}, los \emph{árboles de decisión de profundidad media constante},\footnote{\emph{Read-once decision trees} y \emph{decision trees of constant average depth} respectivamente, en la literatura en inglés.} las funciones aleatorias, etc. (ver las referencias incluidas en \cite{BS25}). Recordemos que para funciones monótonas la conjetura continúa abierta. Por último, animamos al lector interesado a consultar los excelentes artículos sobre las aplicaciones de la conjetura EFI y los últimos progresos en \cite{WWW14}, \cite{KKL+20} y~\cite{Kal20}.

\section*{Apéndice: Aline Bonami, apuntes biográficos}

Aline Bonami es una reconocida matemática francesa experta en varias áreas del análisis matemático. Estudió en la Universidad Paris-Sud, Orsay, siendo la primera estudiante de doctorado de Yves Meyer,\footnote{Yves Meyer recibió entre otros el Premio Princesa de Asturias en 2020 y el Premio Abel en el 2017.} con una tesis titulada \emph{Étude des coefficients de Fourier des fonctions de $L^p(G)$}, donde investiga los multiplicadores de Fourier y las desigualdades de hipercontractividad.
Fue profesora de la Universidad de Orleans desde 1973 hasta su jubilación en el 2006, y desde entonces es profesora emérita.

Durante los primeros años de su carrera, Aline Bonami fue investigadora a tiempo completo del \emph{Centre National de la Recherche Scientifique} (CNRS) en Orsay. Allí se benefició de un ambiente intelectual excepcional en el grupo de análisis armónico liderado por Jean-Pierre Kahane, y tuvo la oportunidad de interactuar con visitantes extranjeros. En particular, fue profundamente influenciada por los cursos dictados por Elias Stein, Raphy Coifman, y Guido Weiss.

Sus trabajos más notables involucran las desigualdades de hipercontractividad, los procesos brownianos fraccionarios, los operadores de Hankel, las proyecciones de Bergman y de Szeg\H{o} y el principio de incertidumbre. Pero también ha hecho incursiones en ámbitos más aplicados como la teoría de señales, la teoría de ondículas y el diagnóstico precoz de la osteoporosis.

\begin{wrapfigure}{l}{0.5\textwidth}
\centering
\includegraphics[width=0.5\textwidth]{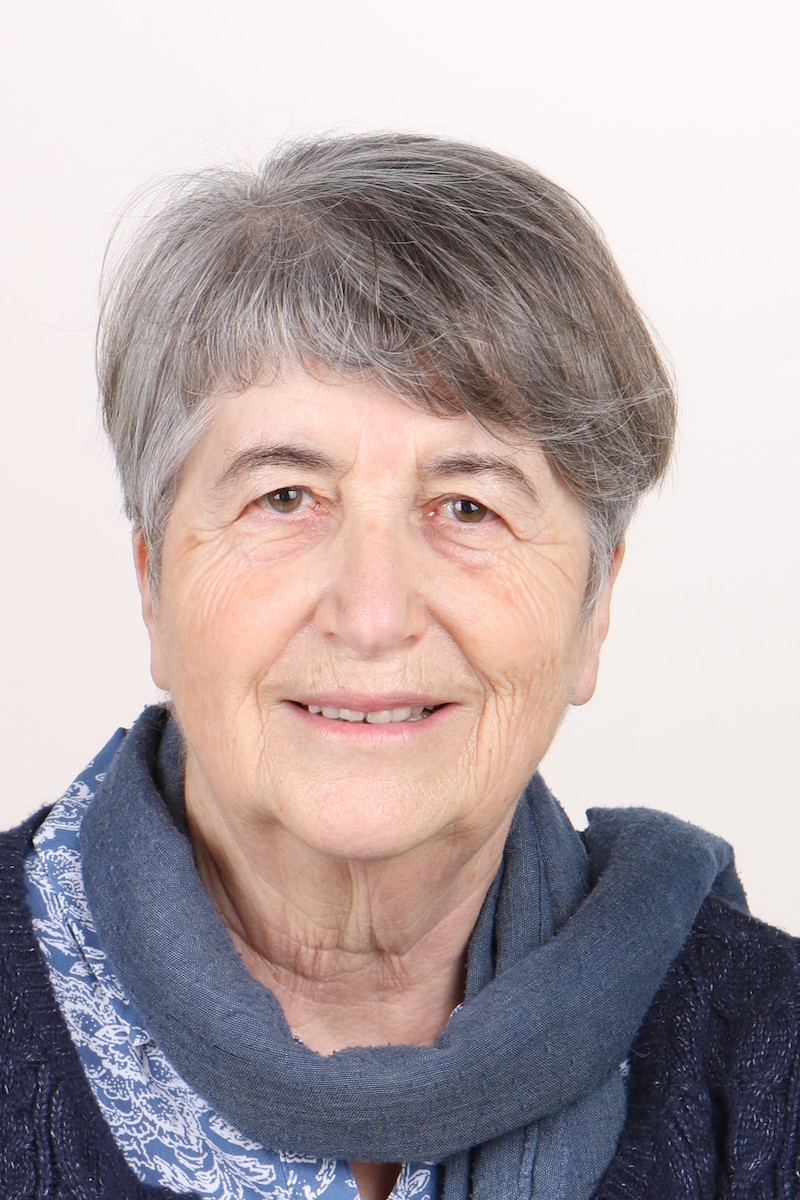}
\pseudocaption{Aline Bonami en 2021. Foto de Aline Bonami.}
\end{wrapfigure}
Aline Bonami tiene una ilustre hoja de servicios. Durante los años 2003-2006 fue la directora científica de matemáticas a cargo de la evaluación en el Ministerio de Investigaciones francés. Fue presidenta de la Sociedad Matemática Francesa del 2012 al 2013. Además ha recibido multiples honores nacionales e internacionales. En el año 2001 recibió el premio Petit D'Ormoy, Carrière, Thébault de la Academia Francesa de Ciencias por sus resultados sobre las proyecciones de Bergman y de Szeg\H{o}, sobre los operadores de Hankel en varias variables complejas y sobre las desigualdades de hipercontractividad. En el año 2002, la Universidad de Gotemburgo, en Suecia, le concedió un doctorado \emph{honoris causa}.
En el año 2005, recibió el Premio del Ministerio Polaco de Educación Nacional para la Investigación y Colaboración, y ese mismo año recibió la \emph{Ordre des Palmes Académiques} (establecida por Napoleón en 1808) en el grado más alto de Comandante. Cinco años más tarde, en el 2010, fue nombrada \emph{Dama de la Legión de Honor} (distinción también creada por Napoleón en 1802). En el año 2014 se organizó en Orleans una conferencia en su honor. En el año 2020 recibió el Premio Stefan Bergman,\footnote{Esa edición de 2020 fue la última del Premio Bergman. En 2023 la American Mathematical Society creó en su lugar la Beca Stefan Bergman.} otorgado por la American Mathematical Society (AMS), junto con Peter Ebenfelt.

La motivación por el premio dice que \emph{<<Aline Bonami recibe el premio por sus contribuciones altamente influyentes en varias variables complejas y en los espacios analíticos. En particular, es reconocida por su trabajo fundamental en las proyecciones de Bergman y de Szeg\H{o} y sus correspondientes espacios de funciones holomorfas. El trabajo de Bonami ha tenido un impacto a largo plazo en la teoría de varias variables complejas, en la teoría de operadores, y en el análisis armónico, y sigue teniendo una gran influencia en la investigación actual que se realiza en estas áreas>>}.

Aline Bonami ha sido muy activa en el ámbito internacional educativo, organizando escuelas CIMPA en Argentina y en Camerún, participando en el comité científico de la Fundación Simons para sus programas en África, así como en el comité científico de los Premios Europeos para Jóvenes Científicos, entre muchas otras actividades.

Para más información sobre la trayectoria de Aline Bonami, se puede consultar:

\begin{itemize}
\item \emph{Mathematics People. Bonami and Ebenfelt awarded 2020 Bergman Prizes.} Notices Amer. Math. Soc. \textbf{68} (2021) (4), 648--649,
\url{https://www.ams.org/journals/notices/202104/rnoti-p648.pdf}

\item Aline Bonami, Wikipedia. \url{https://fr.wikipedia.org/wiki/Aline_Bonami}.

\item A. Bonami, \emph{Why I became a mathematician?}, en \emph{I, mathematician}, P. Casazza, S. G. Krantz y R. D. Ruden (eds.), Mathematical Association of America, Washington, 2015.

\item \url{https://femmes-et-maths.fr/femmes-en-maths/femmes-en-maths/}.
\end{itemize}

\subsubsection*{Agradecimientos}

María José González ha sido financiada por el Ministerio de Ciencia e Innovación (referencia PID2021-123151NB-I00), y por el proyecto \emph{Teoría de Operadores: una aproximación interdisciplinar} (referencia ProyExcel\_00780), proyecto financiado en la convocatoria 2021 de Ayudas a Proyectos de Excelencia, en régimen de concurrencia competitiva, destinadas a entidades calificadas como Agentes del Sistema Andaluz del Conocimiento, en el ámbito del Plan Andaluz de Investigación, Desarrollo e Innovación (PAIDI 2020). Consejería de Universidad, Investigación e Innovación de la Junta de Andalucía.



\phantomsection

\end{document}